\documentclass[10pt, reqno]{amsart}

\usepackage{amssymb}
\usepackage{hyperref} 
\usepackage{mathrsfs,bbold}
\usepackage{amsmath}
\usepackage{amsthm}
\usepackage{enumerate}
\usepackage{dsfont}
\usepackage{color}
\usepackage[a4paper]{geometry}
\usepackage[utf8]{inputenc}
\usepackage{stmaryrd}
\usepackage{titlesec}
\usepackage{mathabx}
\usepackage{appendix}
\usepackage{todonotes, fancyhdr}

\usepackage{tikz-cd}
\usetikzlibrary{calc}

\usepackage[all]{xy}
\let\oldtocsection=\tocsection
\let\oldtocsubsection=\tocsubsection
\let\oldtocsubsubsection=\tocsubsubsection

\renewcommand{\tocsection}[2]{\hspace{0em}\oldtocsection{#1}{#2}}
\renewcommand{\tocsubsection}[2]{\hspace{1em}\oldtocsubsection{#1}{#2}}
\renewcommand{\tocsubsubsection}[2]{\hspace{2em}\oldtocsubsubsection{#1}{#2}}

\numberwithin{theorem}{section}
\numberwithin{equation}{section}

\newcommand{\ssk}{\smallskip}

\renewcommand{\epsilon}{\varepsilon}

\newcommand\bbD{\mathbb{D}}
\newcommand\bbE{\mathbb{E}}

\newcommand\bbN{\mathbb{N}}
\newcommand\bbP{\mathbb{P}}
\newcommand\bbR{\mathbb{R}}

\newcommand\bbZ{\mathbb{Z}}

\newcommand{\mcA}{\mathcal{A}}
\newcommand{\mcB}{\mathcal{B}} 
 
\newcommand{\mcD}{\mathcal{D}}
\newcommand{\mcE}{\mathcal{E}}
\newcommand{\mcF}{\mathcal{F}}
\newcommand{\mcG}{\mathcal{G}}

\newcommand{\mcI}{\mathcal{I}}
\newcommand{\mcJ}{\mathcal{J}}

\newcommand\mcS{\mathcal{S}}

\newcommand{\mcR}{\mathcal{R}}

\newcommand{\bfA}{\mathbf{A}}

\newcommand{\bfG}{\mathbf{G}}

\newcommand{\bfT}{\mathbf{T}}

\newcommand{\scT}{\ensuremath{\mathscr{T}}}
\newcommand{\scW}{\ensuremath{\mathscr{W}}}

\newcommand{\bfa}{\mathbf{a}}
\newcommand{\bfb}{\mathbf{b}}
\newcommand{\bfc}{\mathbf{c}}

\newcommand{\spa}{\mathop{\text{\rm span}}}

\newcommand{\id}{\mathop{\text{\rm id}}}

\newcommand{\tri}{|\!|\!|}

\newcommand{\rest}[1]{\llparenthesis\, #1 \,\rrparenthesis}

\titleformat{\section}[block]
{\filcenter\normalfont\sffamily\bfseries\Large}
{{\hspace{-0.7cm}}\thesection \hspace{0.2em} --\vspace{0.15cm}}{0.5em}{}
\titleformat{\subsection}[runin]
{\filcenter\normalfont\sffamily\bfseries\large}  						  % Without \filcenter
{\hspace{0cm}\thesubsection \hspace{0.5em}--\vspace{0cm}}{.5em}{}  %\hspace{-1cm}
\titlespacing{\subsection}{-0pc}{1.5ex plus .1ex minus .2ex}{0pc}

\titleformat{\subsubsection}[runin]
{\filcenter\normalfont\bfseries\sffamily}
{{\hspace{-0cm}}{\thesubsubsection} \hspace{-0.1em} \hspace{-0.2cm}}{0.5em}{}  

\numberwithin{subsection}{section}
\numberwithin{subsubsection}{subsection}

\newtheoremstyle{mystyle}
{3pt}               %space above
{3pt}               %space below
{\it }                      %bodyfont
{}                      %indent
{\sffamily\bfseries}             %headfont
{}                      %punctuation
{0.5em}                 %space after head
{{#2. }#1 { --}}   % {\llap{#2. }#1{$\;$ --}}

\theoremstyle{mystyle}

\newtheorem{thm}{Theorem}
\newtheorem{lem}[thm]{\hspace{-0.2cm}  {Lemma} }%[section]
\newtheorem{defn}[thm]{ \hspace{-0.3cm} {Definition}}%[chapter]
\newtheorem*{defn*} {Definition}
\newtheorem*{prop*} {Proposition}
\newtheorem*{lem*} {Lemma}
\newtheorem*{cor*} {Corollary}

\newtheoremstyle{mystyle2}
{3pt}               %space above
{3pt}               %space below
{\it }                      %bodyfont
{}                      %indent
{\sffamily\bfseries}             %headfont
{}                      %punctuation
{0.5em}                 %space after head
{\llap{#2 }#1{\hspace{0.2cm}--}}
\theoremstyle{mystyle2}

\newtheorem*{definition*}{Definition}
\newtheorem*{lemma*}{Lemma}
\newtheorem*{thm*}{Theorem}
\numberwithin{equation}{section} % Une \UTF{00E9}quation p.q. est la q \`eme \UTF{00E9}quation de la section p.
\newtheorem*{Remark*}{Remark}

\newcommand\restr[2]{{% we make the whole thing an ordinary symbol
  \left.\kern-\nulldelimiterspace % automatically resize the bar with \right
  #1 % the function
  \vphantom{\big|} % pretend it's a little taller at normal size
  \right|_{#2} % this is the delimiter
  }}

\makeatletter
\newcommand*{\defeq}{\mathrel{\rlap{%
                     \raisebox{0.25ex}{$\m@th\cdot$}}%
                     \raisebox{-0.25ex}{$\m@th\cdot$}}%
                     =}
\makeatother

\makeatletter
\newcommand*{\eqdef}{=\mathrel{\rlap{%
                     \raisebox{0.25ex}{$\m@th\cdot$}}%
                     \raisebox{-0.25ex}{$\m@th\cdot$}}%
                     }
\makeatother

\usepackage{tikz-cd}
\usetikzlibrary{calc}
\usetikzlibrary{shapes.misc}
\usetikzlibrary{shapes.symbols}
\usetikzlibrary{snakes}
\usetikzlibrary{decorations}
\usetikzlibrary{decorations.markings}

\makeatletter
\pgfdeclareshape{crosscircle}
{
  \inheritsavedanchors[from=circle] % this is nearly a circle
  \inheritanchorborder[from=circle]
  \inheritanchor[from=circle]{north}
  \inheritanchor[from=circle]{north west}
  \inheritanchor[from=circle]{north east}
  \inheritanchor[from=circle]{center}
  \inheritanchor[from=circle]{west}
  \inheritanchor[from=circle]{east}
  \inheritanchor[from=circle]{mid}
  \inheritanchor[from=circle]{mid west}
  \inheritanchor[from=circle]{mid east}
  \inheritanchor[from=circle]{base}
  \inheritanchor[from=circle]{base west}
  \inheritanchor[from=circle]{base east}
  \inheritanchor[from=circle]{south}
  \inheritanchor[from=circle]{south west}
  \inheritanchor[from=circle]{south east}
  \inheritbackgroundpath[from=circle]
  \foregroundpath{
    \centerpoint%
    \pgf@xc=\pgf@x%
    \pgf@yc=\pgf@y%
    \pgfutil@tempdima=\radius%
    \pgfmathsetlength{\pgf@xb}{\pgfkeysvalueof{/pgf/outer xsep}}%  
    \pgfmathsetlength{\pgf@yb}{\pgfkeysvalueof{/pgf/outer ysep}}%  
    \ifdim\pgf@xb<\pgf@yb%
      \advance\pgfutil@tempdima by-\pgf@yb%
    \else%
      \advance\pgfutil@tempdima by-\pgf@xb%
    \fi%
    \pgfpathmoveto{\pgfpointadd{\pgfqpoint{\pgf@xc}{\pgf@yc}}{\pgfqpoint{-0.707107\pgfutil@tempdima}{-0.707107\pgfutil@tempdima}}}
    \pgfpathlineto{\pgfpointadd{\pgfqpoint{\pgf@xc}{\pgf@yc}}{\pgfqpoint{0.707107\pgfutil@tempdima}{0.707107\pgfutil@tempdima}}}
    \pgfpathmoveto{\pgfpointadd{\pgfqpoint{\pgf@xc}{\pgf@yc}}{\pgfqpoint{-0.707107\pgfutil@tempdima}{0.707107\pgfutil@tempdima}}}
    \pgfpathlineto{\pgfpointadd{\pgfqpoint{\pgf@xc}{\pgf@yc}}{\pgfqpoint{0.707107\pgfutil@tempdima}{-0.707107\pgfutil@tempdima}}}
  }
}
\makeatother

\colorlet{symbols}{black}    
\colorlet{testcolor}{green!60!black}
\colorlet{supcolor}{red!60!black}

\tikzset{
	root/.style={circle, fill=testcolor!70, draw=testcolor, inner sep=1pt, minimum size=0.5mm},
	dot/.style={circle, draw=black, fill=black, inner sep=0pt, minimum size=1mm},
	noise/.style={circle, draw=black, fill=white, inner sep=0pt, minimum size=1mm},
	h/.style={circle, draw=black, fill=black, inner sep=0pt, minimum size=0.3mm},	
	var/.style={circle, fill=white, draw=purple, fill=purple, inner sep=0pt, minimum size=0.6mm},
	bdot/.style={circle, draw=black, fill=white, inner sep=0pt, minimum size=1mm},
	bluedot/.style={circle,draw=blue, fill=blue, inner sep=0pt, minimum size=2mm},
	noiseblue/.style={circle, fill=blue!20, draw=blue, inner sep=0pt, minimum size=1mm},
	dtestfcn/.style={ultra thick, densely dashed, testcolor,shorten >=1pt,shorten <=1pt,<-},
	testfcnx/.style={semithick,testcolor,shorten >=1pt,shorten <=1pt,<-, postaction={decorate,decoration={markings,mark=at position 0.6 with {\drawx}}}},
	testfcn/.style={semithick, testcolor, shorten >=1pt,shorten <=1pt,-},	
	K/.style= {semithick, shorten >=0pt,shorten <=0pt,-},%{semithick,densely dashed,shorten >=1pt,shorten <=1pt,->},
	kdashed/.style= {semithick,densely dashed,shorten >=1pt,shorten <=1pt,<->},
	DK/.style={thick, densely dotted, shorten >=0pt,shorten <=0pt},   % 	DK/.style={semithick, densely dashed, shorten >=1pt,shorten <=1pt},
	ksup/.style={thick, supcolor, shorten >=1pt,shorten <=1pt},
	dots/.style={semithick,dotted,shorten >=1pt,shorten <=1pt},
	Deps/.style={semithick,draw=black!25,fill=black!25,shorten >=1pt,shorten <=1pt,->},
	kbase/.style={semithick,dotted,shorten >=1pt,shorten <=1pt,->},
	multx/.style={shorten >=1pt,shorten <=1pt,
		postaction={decorate,decoration={markings,mark=at position 0.5 with {\drawx}}}},
	kernelx/.style={semithick,shorten >=1pt,shorten <=1pt,->,
		postaction={decorate,decoration={markings,mark=at position 0.4 with {\drawx}}}},
	kernel1/.style={->,semithick,shorten >=1pt,shorten <=1pt,postaction={decorate,decoration={markings,mark=at position 0.45 with {\draw[--] (0,-0.1) -- (0,0.1);}}}},
	kernel2/.style={->,semithick,shorten >=1pt,shorten <=1pt,postaction={decorate,decoration={markings,mark=at position 0.45 with {\draw[--] (0.05,-0.1) -- (0.05,0.1);\draw[--] (-0.05,-0.1) -- (-0.05,0.1);}}}},
	kernelBig/.style={semithick,shorten >=1pt,shorten <=1pt,decorate, decoration={zigzag,amplitude=1.5pt,segment length = 3pt,pre length=2pt,post length=2pt}},
	rho/.style={dotted,semithick,shorten >=1pt,shorten <=1pt},
	renorm/.style={shape=circle,fill=white,inner sep=1pt},
	labl/.style={shape=rectangle,fill=white,inner sep=1pt},
	xi/.style={circle,fill=symbols!10,draw=symbols,inner sep=0pt,minimum size=1.2mm},
	xix/.style={crosscircle,fill=symbols!10,draw=symbols,inner sep=0pt,minimum size=1.2mm},
	xib/.style={circle,fill=symbols!10,draw=symbols,inner sep=0pt,minimum size=1.6mm},
	xibx/.style={crosscircle,fill=symbols!10,draw=symbols,inner sep=0pt,minimum size=1.6mm},
	not/.style={circle,fill=symbols,draw=symbols,inner sep=0pt,minimum size=0.5mm},
	>=stealth,
	graydot/.style={circle,fill=gray,inner sep=0pt, minimum size=1mm},
	zero/.style={circle,inner sep=0pt, minimum size=1mm, draw},
	kernelprimeeps/.style={densely dashed, semithick,shorten >=1pt,shorten <=1pt},
	smalldot/.style={circle,fill=black,draw=black, solid,inner sep=0pt,minimum size=0.5mm},
	}

\begin{document}

\begin{center}
{\LARGE\sffamily{\textbf{Renormalization of random models: a review}   \vspace{0.5cm}}}
\end{center}

\begin{center}
{\sf I. BAILLEUL\footnote{Univ Brest, CNRS, LMBA - UMR 6205, F- 29238 Brest, France. {\it E-mail}: ismael.bailleul@univ-brest. Partial support from the ANR via the ANR-22-CE40-0017 grant.} 
and
M. HOSHINO\footnote{Graduate School of Engineering Science, Osaka University, Japan. {\it E-mail}: hoshino@sigmath.es.osaka-u.ac.jp. Partial support from the JSPS KAKENHI Grant Number 23K12987.} }
\end{center}

\vspace{1cm}

\begin{center}
\begin{minipage}{0.8\textwidth}
\renewcommand\baselinestretch{0.7} \scriptsize \textbf{\textsf{\noindent Abstract.}} We give a review of three works on the construction of random models for singular stochastic partial differential equations within the theory of regularity structures.
\end{minipage}
\end{center}

%\bigskip

%-----------------------%
\section{Introduction}
\label{SectionIntro}
%-----------------------%

The setting of regularity structures allows to decouple the probabilistic features of the study of a (system of) singular stochastic partial differential equation (SPDE) from its analytic features. The naive formulation of such an equation involves some classical functions or distributions but fails to make sense due to some ill-defined operations involved in the equation: typically a product of two distributions. Think of the archetypal example given by the two-dimensional parabolic Anderson model equation set on the torus
$$
(\partial_t-\Delta) u = u\xi,
$$
for a periodic space white noise $\xi$. The latter has Besov-H\"older regularity $-1-\eta$ for any $\eta>0$, and no better, which gives $u$ some a priori regularity no better than $1-\eta$. This is not sufficient for making sense of the product $u\xi$ as a continuous function of $u$ and $\xi$. One then needs a non-naive way of making sense of the equation even before asking whether it has a (unique) solution or not. The regularity structure formulation of a given singular SPDE is a family of equations
$$
\boldsymbol{u} = F^{\sf M}(\boldsymbol{u})
$$
involving some non-classical objects $\boldsymbol{u}$ called {\it modelled distributions}, and indexed by some parameter $\sf M$ called a {\it model}. Under some mild assumptions on the initial singular equation, each equation of the model-dependent family of equations has a unique `local' solution in a model-dependent space of modelled distributions. To make a link with the initial stochastic equation one would like to apply this local well-posedness result for a random model satisfying some natural constraint involving the random noise in the equation. We could then declare that the associated local solution is a solution to the singular SPDE. Building a random variable taking its values in the space of models and satisfying the above mentioned constraint turns out to be fairly elaborate. This is what renormalization is about.

The analytic and algebraic foundations of regularity structures were laid in the foundational works \cite{Hairer} by M. Hairer  and \cite{BHZ} by Bruned, Hairer \& Zambotti. The first construction of a random model satisfying the above mentioned constraint for a very large class of equation was given by Chandra \& Hairer in \cite{CH}. Its dynamical interpretation was explained in Bruned, Chandra, Chevyrev \& Hairer's work \cite{BCCH}.

The work \cite{CH} holds under some general assumptions on the noise that involve its cumulants, but this work is difficult and technically demanding. It is then fortunate that a different approach to singular SPDEs as a whole, and to renormalization in particular, was developed recently by Otto and co. \cite{OSSW, LOT, LOTT, Tempelmayr} with in mind the study of a number of quasilinear equations. Although the architecture of the approach developed in these works is similar to the architecture of regularity structures, the details differ in a substantial way. The algebra is in particular different, as expected from the point of view adopted, which somehow considers as single objects some sums of objects that are considered separately in a regularity structures setting. Interestingly for us here, they introduce in \cite{LOTT} a very different way of constructing their equivalent of the random model constructed by Chandra \& Hairer in \cite{CH}. By implementing a clever and intricate inductive mechanics, they are able to renormalize their models without renormalizing each of its `pieces' as in the approach of \cite{CH}. The joint consideration of the model together with its Malliavin derivative plays a crucial role in their construction when coupled with the assumption that the law of the random noise satisfies a Poincar\'e-type spectral gap inequality. This type of assumption is somewhat orthogonal to the assumptions of \cite{CH}, and none implies the other.

The clean setting of regularity structures enabled Hairer \& Steele \cite{HS} to implement the mechanics of \cite{LOTT} in a classical regularity structure setting and reprove that the BPHZ renormalization scheme of \cite{BHZ} allows to construct a random model as in \cite{CH}, under the assumption that the noise satisfies a spectral gap inequality. This transfer is not free, though, and Hairer \& Steele introduced for their purpose a notion of pointed modelled distribution that plays a crucial technical role in their construction. This setting is not the only one possible. Bailleul \& Hoshino introduced in \cite{RandomModel} a different way of constructing inductively a large family of random models using a notion of regularity-integrability structure. This construction contains the construction of \cite{HS} as a particular case.

\ssk

The present work aims at giving an overview of the works \cite{CH, HS, RandomModel} by emphasizing some of their mechanics, differences and similarities. We assume the reader is acquainted with the basics of regularity structures so we only recall a number of elementary definitions without going deep in their motivation. The reader will find in \cite{CorwinShen, RSGuide, Berglund} some introductions to the analytic and algebraic aspects of regularity structures. The expository works \cite{HairerBPHZ, HairerTakagi} of Hairer deal specifically with the renormalization of models and some related subjects. We choose to leave aside the very important work \cite{LOTT} of Linares, Otto, Tempelmayr \& Tsatsoulis and concentrate here on the works \cite{CH, HS, RandomModel} as they share the same technical background. For a reader interested in the multi-indices approach to regularity structures we recommend the lecture notes \cite{LO, BOT}.

\medskip

{\bf Organization of this work.} We use Section \ref{SectionBasics} to recall some basic facts and definitions about regularity structures and set a number of notations. A reader acquainted with the subject can directly go to the next section. Section \ref{SectionCH} deals with the diagramatic work \cite{CH} and emphasizes the mechanics at the heart of the BPHZ algorithm for `taming the infinities' that are involved in the construction of the so-called BPHZ model. Section \ref{SectionNondiagramatic} deals with the non-diagramatic works \cite{HS, RandomModel}. We describe the work \cite{HS} of Hairer \& Steele in Section \ref{SectionHS}, where we explain in particular the role played by the pointed modelled distributions. Last we describe in Section \ref{SectionBH} the general construction of \cite{RandomModel} and the flexible setting provided by the notion of regularity-integrability structures.

\bigskip

%-----------------------------------------------------------------------------------------------------%
\section{Basics on regularity structures and a convergence theorem for models}
\label{SectionBasics}
%-----------------------------------------------------------------------------------------------------%

This section is the occasion to introduce a number of notations that will be used at several places below. We set the stage to formulate a version of the convergence theorem for the BPHZ random models that were first introduced by Bruned, Hairer \& Zambotti in \cite{BHZ}.

%%----------------------------------------------%%
\subsection{Basics on regularity structures$\boldmath{.}$ \hspace{0.1cm}}
\label{SubsectionBasicsRS}
%%----------------------------------------------%%

Fix an integer dimension $d\geq 1$. To avoid the technical complexities associated with some anisotropic scaling we consider the isotropic, Euclidean, norm 
$$
\|x\| \defeq \bigg(\sum_{i=1}^d|x_i|^2\bigg)^{\frac12}.
$$ 
The degree of a multiindex $k\in\bbN^d$ is define by 
$$
|k| \defeq \sum_{i=1}^d|k_i|.
$$ 
A {\bf regularity structure} $\scT = \big( \bfA,\bfT,\bfG \big)$ consists of
\begin{itemize} \setlength{\itemsep}{0.1cm}
\item[(1)]
$\bfA$: a subset of $\bbR$ such that the set $\{\alpha\in \bfA\, ;\, \alpha<\gamma\}$ is finite for every $\gamma\in\bbR$.
\item[(2)]
$\bfT=\bigoplus_{\alpha\in \bfA}\bfT_\alpha$: an algebraic sum of Banach spaces $(\bfT_\alpha,\|\cdot\|_\alpha)$.
\item[(3)]
$\bfG$: a group of continuous linear operators on $\bfT$ such that, for any $\Gamma\in \bfG$ and $\alpha\in\bfA$,
$$
(\Gamma-\id)\bfT_\alpha\subset \bfT_{<\alpha}\defeq\bigoplus_{\beta\in \bfA,\, \beta<\alpha}\bfT_\beta.
$$
\end{itemize}
The smallest element $\alpha_0$ of $\bf A$ is called the \emph{regularity} of $\scT$. For any $\alpha\in\bfA$ we denote by $P_\alpha:\bfT\to\bfT_\alpha$ the canonical projection and abuse notation writing for any $\tau\in \bfT$
$$
\|\tau\|_\alpha\defeq\|P_\alpha\tau\|_\alpha.
$$
For any positive integer $r$ we denote by $\mcB_r$ the set of smooth functions $\varphi:\bbR^d\to\bbR$ supported in the unit ball centered at $0$ and such that $\sum_{|k|\le r}\|\partial^r\varphi\|_{L^\infty}\le1$. For $\varphi\in\mcB_r$, $x\in\bbR^d$ and $\lambda\in(0,1]$ we define the function $\varphi_x^\lambda:\bbR^d\to\bbR$ by 
$$
\varphi_x^\lambda(y)  \defeq \lambda^{-d}\varphi\big(\lambda^{-1}(y-x)\big).
$$
(We would use a different definition in a non-isotropic situation.) It converges as $\lambda$ goes to $0$ to a Dirac mass at $x$. The role played by monomials in a classical setting is played here by what we call a model.

\medskip

\begin{defn} \label{HS:def:model}
Let $\scT=(\bfA,\bfT,\bfG)$ be a regularity structure of regularity $\alpha_0$. Fix a positive integer $r>|\alpha_0|$. A {\bf model} ${\sf M}=(\Pi,\Gamma)$ for $\scT$ on $\bbR^d$ consists of two families of continuous linear operators
$$
\Pi = \big\{ \Pi_x : \bfT\to \mcD'(\bbR^d) \big\}_{x\in\bbR^d},\qquad
\Gamma = \{\Gamma_{yx}\}_{x,y\in\bbR^d}\subset \bfG
$$
satisfying the following properties.
\begin{itemize} \setlength{\itemsep}{0.1cm}
\item[(1)]
One has $\Pi_x\Gamma_{xy}=\Pi_y$, $\Gamma_{xx}=\id$ and $\Gamma_{xy}\Gamma_{yz}=\Gamma_{xz}$, for any $x,y,z\in\bbR^d$.
\item[(2)]
For any $\gamma>0$ and any compact subset $C$ of $\bbR^d$ one has 
\begin{align}\label{HS:eq:modelnorm}
\begin{aligned}
\|\Pi\|_{\gamma;C}&\defeq
\max_{\alpha\in\bfA,\,\alpha<\gamma}\sup_{\varphi\in\mcB_r}\sup_{\lambda\in(0,1]}\sup_{x\in C}\sup_{\tau\in\bfT_\alpha\setminus\{0\}}
 \lambda^{-\alpha}
\frac{|(\Pi_x\tau)(\varphi_x^\lambda)|}{\|\tau\|_\alpha}<\infty,
\\
\|\Gamma\|_{\gamma;C}&\defeq
\max_{\alpha,\beta\in\bfA,\,\beta<\alpha<\gamma}
\sup_{x,x+y\in C,\,y\neq0}\sup_{\tau\in\bfT_\alpha\setminus\{0\}}
\frac{\|\Gamma_{(x+y)x}(\tau)\|_\beta}
{\|y\|^{\alpha-\beta}\|\tau\|_\alpha}<\infty.
\end{aligned}
\end{align}
\end{itemize}
We write 
$$
\tri {\sf M}\tri_{\gamma;C} \defeq \|\Pi\|_{\gamma;C}+\|\Gamma\|_{\gamma;C}.
$$ 
Moreover for any two models ${\sf M}^{(1)}$ and ${\sf M}^{(2)}$ we define the quasi-metric $\tri {\sf M}^{(1)} \,,\, {\sf M}^{(2)}\tri_{\gamma;C}$ by replacing $\Pi_x$ and $\Gamma_{(x+y)x}$ in \eqref{HS:eq:modelnorm} with $\Pi_x^{(1)}-\Pi_x^{(2)}$ and $\Gamma_{(x+y)x}^{(1)}-\Gamma_{(x+y)x}^{(2)}$ respectively.
\end{defn}

\medskip

The $\Pi_x(\tau)$ provide the basis for some local expansion of  some class of functions/distributions near an arbitrary point $x\in\bbR^d$. The $\Gamma$ map is used to define a notion of regularity based on some local expansion property for a particular class of functions or distributions that comes in the form of the following definition.

\medskip

\begin{defn} \label{HS:def:MD}
Let $\scT=(\bfA,\bfT,\bfG)$ be a regularity structure and let ${\sf M}=(\Pi,\Gamma)$ be a model for $\scT$ on $\bbR^d$. For any $\gamma\in\bbR$ we define the space $\mcD^\gamma=\mcD^\gamma(\Gamma)$ of {\bf $\gamma$-regular modelled distributions} as the space of all functions $f:\bbR^d\to\bfT_{<\gamma}$ which satisfy
\begin{align}
\label{HS:eq:MD1}
\rest{f}_{\gamma;C}&\defeq
\max_{\alpha\in\bfA,\,\alpha<\gamma}\sup_{x\in C}\|f(x)\|_{\alpha}<\infty,
\\
\label{HS:eq:MD2}
\|f\|_{\gamma;C}&\defeq
\max_{\alpha\in\bfA,\,\alpha<\gamma}
\sup_{x,x+y\in C,\,0<\|y\|\le1}
\frac{\big\|f(x+y) - \Gamma_{(x+y)x}(f(x))\big\|_\alpha}
{\|y\|^{\gamma-\alpha}}<\infty
\end{align}
for any compact subset $C\subset\bbR^d$.
We write 
$$
\tri f\tri_{\gamma;C}\defeq\rest{f}_{\gamma;C}+\|f\|_{\gamma;C}.
$$
Moreover for any two models ${\sf M}^{(1)}$ and ${\sf M}^{(2)}$ and $f^{(i)}\in\mcD^\gamma(\Gamma^{(i)})$ with $i\in\{1,2\}$, we define the quasi metric $\tri f^{(1)},f^{(2)}\tri_{\gamma;C}$ by replacing $f(x)$ in \eqref{HS:eq:MD1} with $f^{(1)}(x)-f^{(2)}(x)$ and replacing $f(x+y)-\Gamma_{(x+y)x}(f(x))$ in \eqref{HS:eq:MD2} with $\big\{f^{(1)}(x+y)-\Gamma_{(x+y)x}^{(1)}(f^{(1)}(x))\big\} - \big\{f^{(2)}(x+y)-\Gamma_{(x+y)x}^{(2)}(f^{(2)}(x))\big\}$.
\end{defn}

\medskip

It is a fundamental fact that one can associate to a modelled distribution a unique function/distribution that is locally well approximated by the model object near every point. This is the content of the {\it reconstruction theorem}, which was first proved by Hairer in Theorem 3.10 of \cite{Hairer}.

\medskip

\begin{thm*}
Let $\scT=(\bfA,\bfT,\bfG)$ be a regularity structure and let ${\sf M} = (\Pi,\Gamma)$ be a model for $\scT$ on $\bbR^d$. For any $\gamma>0$, there exists a unique continuous linear map $\mcR^{\sf M}:\mcD^\gamma(\Gamma)\to\mcD'(\bbR^d)$ satisfying for any compact subset $C$ of $\bbR^d$
$$
\sup_{\varphi\in\mcB_r}\sup_{\lambda\in(0,1]}\sup_{x\in C}\lambda^{-\gamma}
\big|\big(\mcR^{\sf M} f-\Pi_xf(x)\big)(\varphi_x^\lambda)\big|\lesssim\|\Pi\|_{\gamma;\overline{C}}\|f\|_{\gamma;\overline{C}},
$$
where 
$$
\overline{C} \defeq \big\{ x+y\,;\, x\in C,\, \|y\|\le 1\big\}.
$$
Moreover the mapping $({\sf M},f)\mapsto\mcR^{\sf M} f$ is locally Lipschitz continuous with respect to the quasi-metrics $\tri {\sf M}^{(1)} \,,\, {\sf M}^{(2)}\tri_{\gamma;\overline{C}}$ and $\tri f^{(1)},f^{(2)}\tri_{\gamma;\overline{C}}$.
\end{thm*}

\bigskip

%%-----------------------------------------------------------------------------%%
\subsection{The archetype convergence result for a BPHZ model$\boldmath{.}$ \hspace{0.1cm}}
\label{SubsectionBasicsBPHZ}
%%-----------------------------------------------------------------------------%%

As in Section 8 of \cite{Hairer}, we consider here for simplicity some regularity structure constructed from a single noise symbol $\Xi$, a single integration symbol $\mcI$, and some Taylor monomials $X_1,\dots,X_d$.
We first define recursively from the following rules a set $\mcF$ of symbols.
\begin{itemize} \setlength{\itemsep}{0.1cm}
\item[--]
The symbols $\Xi,{\bf1},X_1,\dots,X_d$ are elements of $\mcF$.
\item[--]
For any $\tau\in\mcF$ and any $k\in\bbN^d$, the symbol $\mcI_k(\tau)$ is also an element of $\mcF$.
\item[--]
For any $\tau,\sigma\in\mcF$, the symbol $\tau\sigma$ is also an element of $\mcF$, where we postulate the relations
$$
(\tau\sigma)\eta=\tau(\sigma\eta),\qquad
\tau\sigma=\sigma\tau,\qquad
\tau{\bf1}=\tau.
$$
\end{itemize}
The symbol $\mcI$ represents a given integral operator $K$. For simplicity, and to make things concrete in this paragraph, we take for $K$ the Green function $(1-\Delta)^{-1}$ on $\bbR^d$. (More precisely, we have to define $K$ as a truncated Green function as detailed in Lemma 5.5 of \cite{Hairer}. See Section 5 of \cite{Hos23} for a slightly different definition of $K$.) The fact that $K(x,y)=K(x-y)$ is translation invariant plays a crucial role at some point. This recursive definition of some symbols lends itself to a pictorial representation of these symbols via some decorated trees. We denote by $\circ$ the single node tree that represents the noise $\xi$ in the equation. The node $\bullet$ represents the symbol $\bf1$. The multiindex $k=(k_i)_{i=1}^d$ of $X^k\defeq\prod_{i=1}^dX_i^{k_i}$ is represented as a node decoration. For any tree $\tau$, the symbol $\mcI_k(\tau)$ denotes a planted tree given by the grafting of $\tau$ onto a new root via an edge decorated by $k$. The product is represented as a tree product. For example, the symbol $\Xi\mcI_0(X^m\mcI_k(\Xi)\mcI_\ell(\Xi))$ is represented as the following tree.
$$
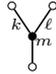
\begin{tikzpicture}[scale=0.3,baseline=0.2cm]
	\node at (0,-0.2)  [noise] (2) {};
	\node at (0,1.15)  [dot] (3) {};
	\node at (0.5,0.9)  {\tiny$m$};
	\node at (-0.9,2.5)  [noise] (left) {};
	\node at (0.9,2.5)  [noise] (right) {};
	\node at (-0.7,1.7) {\tiny $k$};
	\node at (0.7,1.7) {\tiny $\ell$};
	
	\draw[K] (2) to (3);
	\draw[K] (right) to (3);
	\draw[K] (left) to (3);
\end{tikzpicture}
$$
The three nodes $\circ$ and the lower edge have zero decorations, which are omitted.

We fix in this section a parameter $\alpha_0\in(-\infty,-d/2)$ which represents the regularity of the noise in the equation under study. The {\bf degree} map $r:\mcF\to\bbR$ is defined by setting
\begin{equation}\label{Basics:defofr}
\begin{split}
r(X^k)&=|k|,\qquad
r(\Xi)=\alpha_0,\\
r(\mcI_k(\tau))&=r(\tau)+2-|k|,\qquad
r(\tau\sigma)=r(\tau)+r(\sigma).
\end{split}
\end{equation}
The space $\mcF$ contains a number of symbols that are useless for the study of the initial SPDE. The notion of {\bf rule} introduced in Section 5 of \cite{BHZ} is a way to select which symbols will be useful for the regularity structure analysis of the equation. We fix in this section a subset $\mcB\subset\mcF$ consisting of the set of symbols of $\mcF$ which strongly conform to a {\bf complete subcritical rule}. This notion ensures that the set $\{\tau\in\mcB\, ;\, r(\tau)<\gamma\}$ is finite for any $\gamma\in\bbR$, and one can construct a regularity structure $\scT$ by choosing $\bfA=\big\{r(\tau)\,;\,\tau\in\mcB\big\}, \bfT=\spa(\mcB)$ with the direct sum decomposition given by $\bfT_\alpha = \spa\big\{\tau\in\mcB\,;\, r(\tau)=\alpha\big\}$, and choosing for $\bfG$ the set of linear maps $\Gamma:\bfT\to\bfT$ such that
$$
(\Gamma-\id)\tau \in \spa \big\{\sigma\in\mcB\,;\, r(\sigma)<r(\tau)\big\}
$$
for any $\tau\in\mcB$. An {\bf admissible model} is a model for $\scT$ on $\bbR^d$ which interprets the symbol $\mcI$ as the convolution with $K$. The full definition of admissible models is omitted here; it can be found in Definition 8.29 of \cite{Hairer}. Instead, we point only some important properties satisfied by an arbitrary admissible model ${\sf M} = (\Pi,\Gamma)$.
\begin{itemize} \setlength{\itemsep}{0.1cm}
\item[(1)]
There exist some continuous linear maps ${\bf\Pi}:\bfT\to\mcD'(\bbR^d)$ and $F_x\in \bfG$ indexed by $x\in\bbR^d$ such that one has for all $x,y$
$$
\Pi_x = {\bf\Pi}\circ F_x^{-1},\qquad \Gamma_{yx} = F_y\circ F_x^{-1}.
$$
\item[(2)]
The actions of $\Pi_x$ on $X^k$ and $\mcI_k(\tau)$ are given as follows:
$$
\Pi_x(X^k) = (\cdot-x)^k,\qquad
\Pi_x(\mcI_k(\tau)) = K*\Pi_x\tau-\sum_{|\ell|<r(\mcI_k(\tau))} \frac{(\cdot-x)^\ell}{\ell!} \, \big(\partial^\ell K*\Pi_x\tau\big)(x).
$$
\item[(3)]
The maps $F_x$, and consequently $\Gamma_{yx}$, are determined by the map ${\bf\Pi}$.
\end{itemize}

\ssk

%The admissibility condition implies in particular that the map $\Gamma$ is determined by the map ${\bf \Pi}$.
To make things concrete here we will assume that the noise $\xi$ in our equation is the white noise on $\bbR^d$. For any fixed family $\{\varrho^\epsilon\}_{0<\epsilon\leq 1}$ of mollifiers we define the smooth noises 
$$
\xi^\epsilon \defeq \xi*\varrho^\epsilon.
$$
We can then define the family of random admissible models ${\sf M}^\epsilon=(\Pi^\epsilon,\Gamma^\epsilon)$ by the identities
$$
{\bf\Pi}^\epsilon (X^k)(x)=x^k,  \qquad
{\bf\Pi}^\epsilon(\Xi) = \xi^\epsilon,  \qquad
{\bf\Pi}^\epsilon(\mcI_k(\tau)) = \partial^kK*{\bf\Pi}^\epsilon(\tau),  \qquad  
{\bf\Pi}^\epsilon(\tau\sigma) = ({\bf\Pi}^\epsilon\tau)({\bf\Pi}^\epsilon\sigma).
$$
We call it the {\it naive model}. For any tree $\tau$ the smooth function ${\bf \Pi}^\epsilon(\tau)(\cdot)$ is given by an iterated integral with the same structure as the undecorated version of $\tau$. The decorations of $\tau$ inform us about which kernels are used in the iterated integral and which polynomial functions of the integration variables are inserted at any given place. One has for instance
\begin{equation*}
{\bf \Pi}^\epsilon
\bigg( \, 
\begin{tikzpicture}[scale=0.3,baseline=0.2cm]
	\node at (0,-0.2)  [noise] (2) {};
	\node at (0,1.15)  [dot] (3) {};
	\node at (0.5,0.9)  {\tiny$m$};
	\node at (-0.9,2.5)  [noise] (left) {};
	\node at (0.9,2.5)  [noise] (right) {};
	\node at (-0.7,1.7) {\tiny $k$};
	\node at (0.7,1.7) {\tiny $\ell$};
	
	\draw[K] (2) to (3);
	\draw[K] (right) to (3);
	\draw[K] (left) to (3);
\end{tikzpicture}
\, \bigg)(x) = \xi^\epsilon(x) \int K(x,y)\, y^m\, (\partial^k K)(y,a)\xi^\epsilon(a) \, (\partial^\ell K)(y,b)\xi^\epsilon(b)\,da\,db\,dy.
\end{equation*}
%Would the internal node in the above tree $\tau$ have an $X$ decoration we would have 
%$$
%{\bf \Pi}^\epsilon(\tau) = \xi^\epsilon(x) \int K(x,y) y \, (\partial K)(y,a)\xi^\epsilon(a) \, (\partial K)(y,b)\xi^\epsilon(b)\,dadbdy.
%$$
The recentered version $\Pi^\epsilon_x(\tau)$ of ${\bf \Pi}^\epsilon(\tau)$ is of the form 
$$
\Pi^\epsilon_x(\tau) = {\bf \Pi}^\epsilon(\tau) - \sum_\sigma c^\varepsilon(\tau,\sigma)(x)\,{\bf \Pi}^\epsilon(\sigma)
$$
for some subtrees $\sigma$ of $\tau$. It corresponds somehow to replacing a ${\bf \Pi}^\epsilon(\rho)$ by its Taylor remainder at some $\rho$-dependent order and iterating this operation in the tree structure of the iterated integral, from the leaves to the root of the tree, in a multiplicative way. See for instance Equation (21) in \cite{BailleulBrunedLocality} -- the details do not matter here. In any case the quantities ${\bf \Pi}^\epsilon(\tau)$ and $\Pi^\epsilon_z(\tau)$ are some polynomial functionals of the noise $\xi^\epsilon$.  We will see in the introduction of Section \ref{SectionCH} that we cannot expect the convergence of ${\sf M}^\epsilon$ as $\epsilon$ goes to $0$. One has to tweak the naive model to make it converge. A class of continuous linear maps $R:\bfT\to\bfT$ called {\bf preparation maps} was introduced by Bruned for that purpose in \cite{Bru18} -- they were called {\it admissible maps} therein. One can associate to each $0<\epsilon\leq 1$ and each preparation map $R$ an admissible model ${\sf M}^{\epsilon,R}$ as in Section 3 of \cite{Bru18}. The details of this construction are not important here. The so-called BPHZ renormalized model of \cite{BHZ} corresponds to a subclass of preparation maps of the form
\begin{align} \label{HS:eq:preparation}
R_\ell=(\ell\otimes\id)\Delta_r^-,
\end{align}
where $\Delta_r^-$ is a map splitting a given symbol $\tau$ identified with a tree into a subtree with negative degree which contains the root of $\tau$ and the remaining graph, and the map $\ell$ assigns a real number to each tree with negative degree and has some morphism property -- see Section 4.1 of \cite{Bru18} for the definitions. Write ${\sf M}^{\epsilon,\ell},\Pi^{\varepsilon,\ell}$ and ${\bf \Pi}^{\epsilon,\ell}$ instead of ${\sf M}^{\varepsilon,R_{\ell}},\Pi^{\varepsilon,R_\ell}$ and ${\bf \Pi}^{\varepsilon,R_\ell}$ for simplicity. This model satisfies an identity of the form
\begin{equation} \label{EqRenormalizationRelation}
\Pi_x^{\epsilon,\ell}(\tau) = \Pi_x^{\epsilon}(A^\ell\tau) = \Pi_x^\epsilon(\tau) + \cdots
\end{equation}
for all $x$ and $\tau$, for some map $A^\ell$ that is a perturbation of the identity such that $A^\ell-\textrm{id}$ is niplotent. In those terms, Theorem 6.18 in \cite{BHZ} says that for each $0<\epsilon\leq 1$ there exists a unique character $\ell^\epsilon$ such that for any $\tau\in\mcB$ with negative degree and any $x$ we have
$$
\bbE\big[ \big({\bf \Pi}^{\epsilon,\ell^\epsilon}\tau\big)(x)\big] = 0.
$$
The model 
$$ 
\overline{\sf M}^\epsilon \defeq {\sf M}^{\epsilon,\ell^\epsilon}
$$
is called the {\bf BPHZ model}. Its probabilistic convergence as $\epsilon>0$ goes to $0$ takes the following form.

\medskip

\begin{thm} [{\cite[Theorem 6]{BH23}}] \label{HS:thm:main}
Suppose that any $\tau\in\mcB\setminus\{\Xi\}$ satisfies
\begin{align} \label{HS:thm:main:asmp}
r(\tau)>-\frac{d}2.
\end{align}
%{\color{blue} Under some appropriate assumption on the law of the noise,} 
Then for any $\gamma>0, 1\leq q<\infty$ and any compact set $C\subset\bbR^d$ we have
$$
\lim_{\epsilon_1,\epsilon_2\to 0}
\bbE\big[
\tri\overline{\sf M}^{\epsilon_1} \,,\, \overline{\sf M}^{\epsilon_2}\tri_{\gamma;C}^q
\big] = 0.
$$
\end{thm}

\medskip

Here we assume that $\xi$ is the white noise on $\bbR^d$, but the above theorem can be extended to more general noises satisfying some appropriate assumptions on their law. Chandra \& Hairer assume in \cite{CH} some moment bounds on the noise while a spectral gap assumption on the law of the noise is assumed in \cite{HS} and \cite{RandomModel}, following Linares, Otto, Tempelmayr \& Tsatsoulis' work \cite{LOTT}.

\bigskip

%-----------------------------------------------------------------------------------------------------%
\section{Diagramatic method: Chandra \& Hairer \cite{CH} and iterated integrals}
\label{SectionCH}
%-----------------------------------------------------------------------------------------------------%

We explain in this section some of the ingredients of Chandra \& Hairer's work \cite{CH}.

Although this may not be obvious at first sight from the notion of convergence of models given in Definition \ref{HS:def:MD}, in the above white noise setting, and more generally for a large class of Gaussian noises, one can show that a quantitative estimate of the form
\begin{equation} \label{EqEstimateConvergenceModels}
\bigg\Vert \int \Big(\overline{\sf \Pi}^{\epsilon_1}_x\tau - \overline{\sf \Pi}_x^{\epsilon_2}\tau\Big)(y) \, \varphi^\lambda_x(y) \, dy \bigg\Vert_{L^2(\Omega)}^2 \lesssim o_{\max(\epsilon_1,\epsilon_2)}(1)\,\lambda^{2(r(\tau)+\kappa)},
\end{equation}
implies the convergence of the models $\overline{\sf M}^\epsilon$ to some limit model, as stated in Theorem \ref{HS:thm:main}. The notation $o_\rho(1)$ denotes a factor converging to $0$ as $\rho>0$ goes to $0$. We require here that the estimate is uniform in $\varphi\in\mcB_r$, locally uniformly in $x$, for some positive constant $\kappa$, for all the symbols $\tau$ of negative degree $r(\tau)$. One can think of \eqref{EqEstimateConvergenceModels} as the main condition in a form of Kolmogorov regularity theorem in a regularity structure setting -- see Theorem 10.7 in \cite{Hairer}.

The kernels that we typically use in a singular SPDE setting have a logarithmic or polynomial singularity on their diagonal
$$
K(x-y) \simeq \Big(\log\vert x-y\vert \,\; \textrm{or} \,\; \vert x-y\vert^{-\beta}\Big)
$$
for some positive constant $\beta$. This singular feature of the kernels of the operators of interest is the very source of the problem of renormalization. Here is a simple example that shows that the estimate \eqref{EqEstimateConvergenceModels} does not hold for the naive model. For a space white noise $\xi$ on the two-dimensional torus with $\xi^\epsilon$ of covariance $C^{\epsilon}$ converging to a Dirac covariance, and $K(x,y)=\vert x-y\vert^{-\beta}$ for $0<\beta<2$, the expectation of $\big(\hspace{0.02cm}{\sf \Pi}_0^{\epsilon}\,
\begin{tikzpicture}[scale=0.3,baseline=0.06cm]
	\node at (0,0)  [noise] (2) {};
	\node at (0,1.15)  [noise] (3) {};	
	\draw[K] (2) to (3);
\end{tikzpicture} \,
\big)(y)$ is equal to
$$
\bbE\bigg[ \bigg(\int K(y,c)\xi^\epsilon(c)dc - \int K(0,c)\xi^\epsilon(c) \bigg) \xi^\epsilon(y) \bigg] = \int K(y,c) C^\epsilon(c,y) dc - \int K(0,c)\,C^\epsilon(c,y)dc.
$$
The second integral on the two dimensional torus is converging to $K(0,y)=|y|^{-\beta}$ in a distributional sense as $\epsilon>0$ goes to $0$, since $\beta<2$, but the first integral is diverging as $\epsilon^{-\beta}$. So the random variable $\int \big(\overline{\sf \Pi}_x^{\epsilon}\tau\big)(y) \, \varphi^\lambda_x(y) \, dy$ cannot converge in $L^2(\Omega)$; it cannot satisfy a fortiori an estimate of the form \eqref{EqEstimateConvergenceModels}.  It is plain on this example that the source of the problem is in the singular character of the kernel $K$ on its diagonal. Two non-trivial things are thus involved in the estimate \eqref{EqEstimateConvergenceModels}: The fact that it remains $(\epsilon_1,\epsilon_2)$-uniformly finite, and the fact that it behaves as $\lambda^{2(r(\tau)+\kappa)}$ as a function of $\lambda$. We will concentrate here on the first point: taming the infinities.

\ssk

%%--------------------------------------------------%%
\subsection{From models to iterated integrals$\boldmath{.}$ \hspace{0.1cm}}
\label{SubsubsectionGraphs}
%%--------------------------------------------------%%

As a consequence of the translation invariance of the law of the noise, the $L^2(\Omega)$-expectation \eqref{EqEstimateConvergenceModels} is actually independent of $x$. Since we are working with some polynomial functionals of the noise the random variable $\int  \big(\overline{\sf \Pi}^\epsilon_0\tau\big)(a) \, \varphi^\lambda_0(a)\,da$. It has an $L^2(\Omega)$-orthogonal chaos decomposition
\begin{equation} \label{EqChaosDecomposition}
\int  \big(\overline{\sf \Pi}^\epsilon_0\tau\big)(a) \, \varphi^\lambda_0(a)\,da = \sum_{0\leq m\leq n_\Xi(\tau)} I_m\bigg(\int \overline{{\bf W}}^{\epsilon,m}_\tau(a,\cdot)\, \varphi^\lambda_0(a)\,da\bigg)
\end{equation}
encoded by some kernels $\int \overline{{\bf W}}^{\epsilon,m}_\tau(a,\cdot)\, \varphi^\lambda_0(a)\,da$ that are elements of the symmetrized $m$-th tensor product of the $L^2$ space of the state-space. See for instance Section 10.1 of \cite{Hairer}. We wrote here $n_\Xi(\tau)$ for the number of noise symbols in $\tau$. The functions $\overline{{\bf W}}^{\epsilon,m}_\tau(a,\cdot)$ are given by some iterated integrals and, setting 
$$
\overline{{\bf W}}^{(\epsilon_1,\epsilon_2),m}_\tau \defeq \overline{{\bf W}}^{\epsilon_1,m}_\tau - \overline{{\bf W}}^{\epsilon_2,m}_\tau,
$$ 
the squared $L^2(\Omega)$ norm of \eqref{EqChaosDecomposition} is equal to
\begin{equation} \label{EqL2OmegaNormModel}
\overline{(\star)}^{(\epsilon_1,\epsilon_2),m} \defeq \sum_{0\leq m\leq n_\Xi(\tau)} \iiint \overline{{\bf W}}_\tau^{(\epsilon_1,\epsilon_2),m}(a,y) \, \overline{{\bf W}}_\tau^{(\epsilon_1,\epsilon_2),m}(b,y)\, \varphi^\lambda_0(a)\varphi^\lambda_0(b)\,da\,db\,dy
\end{equation}
by It\^o isometry. (In the more general case of a translation invariant Gaussian noise with covariance $C$ we would have an expression of the form $\iint \overline{{\bf W}}(a,y_1) \,C(y_1-y_2)\, \overline{{\bf W}}(b,y_2)dy_1dy_2$ instead of the above integral $\int \overline{{\bf W}}(a,y) \overline{{\bf W}}(b,y)\,dy$.) The convergence criterion \eqref{EqEstimateConvergenceModels} brings us back to proving some quantitative estimates on some iterated integrals.

We denote by ${\bf W}^{\epsilon,m}_\tau$ the equivalent of $\overline{{\bf W}}^{\epsilon,m}_\tau$ for the naive random variable ${\sf \Pi}_0^\varepsilon\tau$.
For example, we consider the tree 
$\tau = \begin{tikzpicture}[scale=0.3,baseline=0.14cm]
	\node at (0,0)  [noise] (2) {};
	\node at (0,0.7)  [dot] (3) {};
	\node at (-0.5,1.4)  [noise] (left) {};
	\node at (0.5,1.4)  [noise] (right) {};
	
	\draw[K] (2) to (3);
	\draw[DK] (right) to (3);
	\draw[DK] (left) to (3);
\end{tikzpicture}$, where the solid line denotes $\mcI_0$ and the dotted line denotes $\mcI_k$ with some $k$.
Then we have
$$
({\bf W}^{\varepsilon,3}_\tau)\big(a,y_1,y_2,y_3\big) = \varrho^\varepsilon(a-y_1)\int K(a,c) (\partial^kK*\varrho^\varepsilon)(c,y_2) (\partial^kK*\varrho^\varepsilon)(c,y_3) \, dc,
$$
${\bf W}^{\varepsilon,2}_\tau=0$, and 
\begin{align*}
({\bf W}^{\varepsilon,1}_\tau)(a,y) &= 2 \int\varrho^\varepsilon(a-b) K(a,c) (\partial^kK*\varrho^\varepsilon)(c,b) (\partial^kK*\varrho^\varepsilon)(c,y) \, db\,dc\\
&\quad+ \varrho^\varepsilon(a-y) \int K(a,c) (\partial^kK*\varrho^\varepsilon)(c,b) (\partial^kK*\varrho^\varepsilon)(c,b) \, db\,dc.
\end{align*}
%Denoting by ${\bf W}^{\epsilon,m}_\tau$ the equivalent of $\overline{{\bf W}}^{\epsilon,m}_\tau$ for ${\sf \Pi}_0\tau$, we have for instance
%$$
%({\bf W}^{0,3}_\tau)\big(a,y_1,y_2,y_3\big) = \delta_a(y_1)\int K(y_1,c) \partial K(c,y_2) \partial K(c,y_3) \, dc,
%$$
%${\bf W}^{0,2}_\tau=\mathcal{W}^0\tau=0$, and 
%$$
%({\bf W}^{0,1}_\tau)(a,y) = 2 \int K(a,c) \partial K(a,c) \partial K(c,y_1) \, dc,
%$$
%on the example of the tree 
%$\tau = \begin{tikzpicture}[scale=0.3,baseline=0.14cm]
%	\node at (0,0)  [noise] (2) {};
%	\node at (0,0.7)  [dot] (3) {};
%	\node at (-0.5,1.4)  [noise] (left) {};
%	\node at (0.5,1.4)  [noise] (right) {};
%	
%	\draw[K] (2) to (3);
%	\draw[DK] (right) to (3);
%	\draw[DK] (left) to (3);
%\end{tikzpicture}$. 
We note from \eqref{EqRenormalizationRelation} that $\overline{{\bf W}}^{\epsilon,m}_\tau$ is a perturbation of ${\bf W}^{\epsilon,m}_\tau$, so we can see \eqref{EqL2OmegaNormModel} as a perturbation of the corresponding quantity $(\star)^{(\epsilon_1,\epsilon_2),m}$ obtained by replacing $\overline{{\bf W}}^{\epsilon,m}_\tau$ by ${\bf W}^{\epsilon,m}_\tau$ therein. This brings us to the question addressed in the next section: {\it Is there a robust recipe for extracting from the a priori diverging quantities $(\star)^{(\epsilon_1,\epsilon_2),m}$ some converging quantities as $\epsilon>0$ goes to $0$?} Ideally, this convergent part of $(\star)^{(\epsilon_1,\epsilon_2),m}$ should be $\overline{(\star)}^{(\epsilon_1,\epsilon_2),m}$. The answer we will provide will give some strong hints as to why the $L^2(\Omega)$ norm \eqref{EqEstimateConvergenceModels} is finite. Obtaining the $\lambda^{2(r(\tau)+\kappa)}$ factor is a different matter on which we will only say a word in Section \ref{SubsectionRecentering}.

\medskip

%%------------------------------------%%
\subsection{The BPHZ mechanics$\boldmath{.}$ \hspace{0.1cm}}
\label{SubsectionBPHZ}
%%------------------------------------%%

It will be useful in the sequel to consider not only some numbers given by some integrals but also some functions given by integrals. We will consider here a class of functions $\bf G$ of $z=(z_1,\dots,z_n)$ indexed by some oriented graphs $G=(\mathcal{E},\mathcal{V})$ with
$$
{\bf G}(z) = {\bf G}(z_1, \dots, z_n) = \int \prod_{i=1}^n\prod_{e\in\mathcal{E}_\partial^i} K_e(z_i-x_{e_-}) \prod_{e'\in\mathcal{E}\backslash\mathcal{E}_\partial} K_{e'}(x_{e'_+}-x_{e'_-})\,dx.
$$
These are our model quantities for $(\star)^{(\epsilon_1,\epsilon_2),m}$, for some fixed values of $\epsilon_1,\epsilon_2$ that do not appear in the notations. We use the notation $dx$ for a shorthand notation for the integration with respect to all the variables that we integrate. Here we denote by $e = (e_+,e_-)\in\mathcal{E}$ a generic edge and write $x_{e_{\pm}}$ for the corresponding variable in the integral. The set $\mathcal{V}$ indexes both the variables $(z_1,\dots,z_n)$ and the integration variables in the integral. To distinguish $(z_1,\dots,z_n)$ from the other variables we talk of the former as the `{\it external variables}'. We denote by $\mathcal{V}_{\textrm{ext}}$ the vertices corresponding to $\big\{z_1,\dots,z_n\big\}$, and we write $\mathcal{V}_{\textrm{int}}$ for the set of `internal' vertices, corresponding to the integration variables. We denote by $\mathcal{E}^i_\partial$ the set of edges for which $x_{e_+}=z_i$, for $1\leq i\leq n$. We set $\mathcal{E}_\partial \defeq \bigsqcup_{i=1}^n \mathcal{E}_\partial^i$ and note that $\prod_{i=1}^n\prod_{e\in\mathcal{E}_\partial^i} = \prod_{e\in\mathcal{E}_\partial}$. The variables associated with the endpoints of an edge in $\mcE\backslash\mcE_\partial$ are thus some internal variables. A positive translation invariant (smooth) kernel $K_e(x-y)$ is associated to each edge $e\in\mathcal{E}$, with diagonal behavior quantified as
\begin{equation} \label{EqControlKernel}
\big\vert \partial^k K_e(z) \big\vert  \lesssim \Vert z\Vert^{-a_e-\vert k\vert}
\end{equation}
for some positive constant $a_e$ and all multiindices $k$ with $\vert k\vert$ bounded by an appropriate constant. In the setting of Section \ref{SubsubsectionGraphs} the smooth kernels $K_e$ implicitly depend on $\epsilon_1,\epsilon_2$ but the above control is assumed to be uniform in $\epsilon_1,\epsilon_2$. We associate to any sub-graph $g=(\mathcal{E}(g),\mathcal{V}(g))$ of $G$ the corresponding function $\bf g$ of its external variables and we define the \textit{superficial degree of divergence} of $g$ as
$$
\omega(g) \defeq d\big( \vert\mathcal{V}_{\textrm{int}}(g)\vert - 1\big) - \sum_{e\in\mathcal{E}(g)} a_e,
$$ 
where $\vert\mathcal{V}_{\textrm{int}}(g)\vert$ is the number of internal vertices of $g$ -- that is the number of variables that one integrates. A graph $g$ for which $\omega(g) \leq 0$ is called \textit{superficially divergent}. We will simply say `{\it divergent}' below. The denomination of $\omega(g)$ is justified by a result of S. Weinberg \cite{Weinberg} from the early sixties that says that if all the subgraphs of $G$ have a positive superficial degree of divergence then the integral that defines formally $\bf G$ is actually absolutely convergent. We set
$$
\omega_+\hspace{-0.03cm}(g) \defeq \max\big(0,-\omega(g)\big).
$$

%%%----------------------------------%%%
\subsubsection{The BPHZ algorithm$\boldmath{.}$}
\label{SubsubsectionBPHZ}
%%%----------------------------------%%%

This procedure is named after Bogoliubov \& Parasiuk for some of their foundational works on the renormalization problem in quantum field theory in the mid fifties, which was clarified by some subsequent works of Hepp and Zimmermann in the late sixties. 

For a sub-graph $g=(\mathcal{E}(g),\mathcal{V}(g))$ of $G=(\mathcal{E}(G), \mathcal{V}(G))$ we denote by $G\backslash g$ the graph obtained from $G$ by removing $\mathcal{E}(g)$ from $\mathcal{E}(G)$. We use the suggestive notation ${\bf G} \backslash {\bf g}$ for the corresponding function and note that the functions $\bf g$ and ${\bf G} \backslash {\bf g}$ have the same external variables; part of them may be some of the variables of $\bf G$. One can write
\begin{equation} \label{EqFactorization}
{\bf G}(z) = \int {\bf g}(z,y) ({\bf G} \backslash {\bf g})(z,y) \, dy,
\end{equation}
keeping in mind that each function $\bf g$ and ${\bf G} \backslash {\bf g}$ may or may not depend on part or all of the variables $z$ of $\bf G$, depending on the situation. We set 
$$
(\mathscr{T}_g{\bf G})=0 \; \textrm{ if } \; \omega(g)>0
$$ 
and otherwise
\begin{equation*}
(\mathscr{T}_g{\bf G})(z) \defeq \int {\bf g}(z,y) \sum_{|k|<\omega_+\hspace{-0.03cm}(g)} \partial_y^k({\bf G} \backslash {\bf g})(z,\widetilde{y}_1)\frac{(y-\widetilde{y}_1)^k}{k!}\,dy,
\end{equation*}
for a point $(y_1,\dots, y_1) \eqdef \widetilde{y}_1$ chosen arbitrarily amongst the external variables of $\bf g$ that we integrate in \eqref{EqFactorization}. We do not emphasize the dependence of the operation $\mathscr{T}_g$ on the somewhat arbitrary choice of point $y_1$ as that choice is of no importance for what follows. We note that if $g$ and $h$ are two disjoint sub-graphs of $G$ then the operators $\mathscr{T}_g$ and $\mathscr{T}_h$ commute. 

\medskip

A family $\mathcal{T}$ of divergent sub-graphs of $G$ such that they are pairwise either disjoint or one is included in the other is called a tree of divergent sub-graphs. The support of such a tree is the union of its sub-graphs; it is a subset of $G$. A forest is a union of trees whose supports are pairwise disjoint. We denote by $\mathscr{F}$ the set of forests of divergent sub-graphs of $G$. The BPHZ {\it renormalization of $\bf G$} is defined as
\begin{equation} \label{EqBPHZ}
\overline{\bf G} \defeq \sum_{\mathcal{F}\in\mathscr{F}} \prod_{g\in \mathcal{F}} (-\mathscr{T}_g){\bf G}.
\end{equation}
For each forest $\mathcal{F}=(\mathcal{T}_1,\dots,\mathcal{T}_\ell)$ the order in the product 
$$
\prod_{g\in \mathcal{F}} = \prod_{1\leq i\leq \ell} \prod_{g\in \mathcal{T}_i}
$$ 
is done within each tree from the leaves to the root. The order in the product over $i$ does not matter since $\mathscr{T}_g$ and $\mathscr{T}_h$ commute if $g$ and $h$ have disjoint supports. Note that for any given forest $\mathcal{F}^*\in\mathscr{F}$ one has
$$
\sum_{\mathcal{F}\in\mathscr{F}, \, \mathcal{F}\subset \mathcal{F}^*} \, \prod_{g\in \mathcal{F}} (-\mathscr{T}_g){\bf G} = \prod_{g\in\mathcal{F}^*}(\textrm{Id}-\mathscr{T}_g){\bf G}.
$$
We expect from the Taylor remainder maps $(\textrm{Id}-\mathscr{T}_g)$ that they produce some terms that make convergent the integral that defines $\bf G$. This idea of using some Taylor remainder maps seems promising, as in Hadamard's finite part extension procedure where we define the distribution $\vert y\vert^{-\beta}$ in a neighborhood of $0$ in $\bbR^d$, for some $\beta>0$, by setting
$$
\big\langle \vert \cdot\vert^{-\beta} \,,\, f\big\rangle \defeq \int \vert y\vert^{-\beta} \bigg\{f(y) - \sum_{\vert k\vert <[\beta-d+1]} (\partial^kf)(0)\,\frac{y^k}{k!}\bigg\} dy,
$$
where $[a]$ stands for the integer part of a real number $a$. We can even think of iterating this kind of operations in some situations where we have some nested divergent subgraphs, `curing' them one after the other in an increasing order for the inclusion relation, as in $\prod_{g\in\mathcal{F}^*}(\textrm{Id}-\mathscr{T}_g)$. But what should we do if two divergent subgraphs are overlapping, with the possibility that one Taylor remainder operation for curing one graph destroys the curing effect that the other operation has on the other graph? It turns out that the BPHZ renormalization prescription \eqref{EqBPHZ} handles this type of situation.

\medskip

%%%-----------------------------------------------------------------------------%%%
\subsubsection{Multiscale analysis and parcimonious renormalization$\boldmath{.}$}
\label{SubsubsectionMultiscale}
%%%-----------------------------------------------------------------------------%%%

A good approach to understand the formula \eqref{EqBPHZ} consists in writing each kernel $K_e(z)$ as a sum of kernels localized in some annuli of dyadic radius $2^{-i}$. We will illustrate this here in the particular case where we work on $\bbR^d$ and all the kernels $K_e$ are equal to the Green function $K$ of the operator $(\Delta-1)$. This kernel has a logarithmic ($d=2$) or polynomial ($d\geq 3$) singularity on its diagonal. We deal with the general case in a very similar way. To stay aligned with the previous picture where we work with some {\it smooth} kernels that are only controlled by \eqref{EqControlKernel} we could work with $K_e = e^{\nu\Delta}(\Delta-1)^{-1}$, for a positive parameter $\nu$ that plays the role of $(\epsilon_1,\epsilon_2)$. We do not loose anything in working directly with $(\Delta-1)^{-1}$, which we do. We thus have
\begin{equation*} \begin{split}
K(z) &= \int_0^\infty \frac{\exp\left(\frac{-|z|^2}{2t}-t\right)}{(4\pi t)^{d/2}}\,dt = \sum_{i\in\bbZ} K^i(z)
\end{split} \end{equation*}
with
\begin{equation} \label{EqBoundKi}
K^i(z) \defeq \int_{2^{-2i}}^{2^{-2(i-1)}} \frac{\exp\left(\frac{-|z|^2}{2t}-t\right)}{(4\pi t)^{d/2}} \, dt \lesssim 2^{(d-2)i}\,\exp\big(-c\,2^{2i}|z|^2\big)
\end{equation}
for some positive constant $c$. We talk of $i$ as the {\it scale} of the kernel $K^i$. The function $\bf G$ can then be represented as a sum
$$
{\bf G} = \sum_\mu {\bf G}^\mu
$$
indexed by the different scale assignments on each edge of the graph $G$. The problem comes from the fact that this sum over the scale assignments has no reason to converge a priori -- uniformly in $\epsilon_1,\epsilon_2$ when we consider the $(\epsilon_1,\epsilon_2)$-dependent situation of Section \ref{SubsubsectionGraphs}.

\ssk

For each scale assignment $\mu$ the problems come from couples of diverging sub-graphs $g_1\subset g_2$, without any other diverging subgraph in between, and such that all the scales of the edges of $g_1$ are greater than the scales of any of the edges in $\mathcal{E}(g_2\backslash g_1)$. One says that $g_1$ {\it is quasi-local for $g_2$ for this scale assignment}. (This notion depends on the scale assignment!) One can indeed prove that the sum of the ${\bf G}^\mu$ over the scale assignments $\mu$ that do not show this `pathology' converges. (Section \ref{SubsubsectionForests} gives some elements about that point.) The following remarks highlight the problem and propose a remedy for it.

\begin{itemize}
	\item[(a)] Denote by $m$ the smallest scale of $g_1$ and by $M$ the largest scale of the edges of $g_2$ that are not some edges of $g_1$. Write
	$$
	{\bf g}_2(z) = \int {\bf g}_1(z,y) ({\bf g}_2 \backslash {\bf g}_1)(z,y) \, dy.
	$$
	The graph $g_1$ is quasi-local for $g_2$ for the scale assignment $\mu$ when $m>M$. In that case  the small graph appears as almost local/pointwise from the point of view of the big graph, in the sense that ${\bf g}_1$ gives from \eqref{EqBoundKi} a $o(1)$ contribution to the integral whenever two of its arguments $y_i, y_j$ satisfy $|y_i-y_j| \gg 2^{-m}$. This is where we loose control on the summation over these bad scale assignments as this basic analysis only gives some control of the integral of order $\sum_{m\geq 0} 2^{-m\omega(g_1)}$ which diverges when $\omega(g_1)\leq 0$.   \vspace{0.1cm}
	
	\item[(b)] For any pair $(g_1,g_2)$ of divergent graphs such that $g_1$ is quasi-local for $g_2$ for the scale assignment $\mu$ one can trade ${\sf g}_2$ for 
	$$
	\big(\textrm{Id} - \mathscr{T}_{g_1}\big){\bf g}_2.
	$$
	With the notations of Section \ref{SubsubsectionBPHZ}, the support of the function ${\bf g}_1(z,y)$ of $y$ is essentially contained in a neighbourhood of the deep diagonal ($y_i=y_j$ for all $i,j$) of size $2^{-m}$. The function $({\bf g}_2\backslash {\bf g}_1)(z,y)$ of $y$ having all its scales smaller than $M$, one gains a factor $2^{-(\omega_+(g_1)+1)(m-M)}$ by replacing the function $({\bf g}_2\backslash {\bf g}_1)(z,\cdot)$ by its Taylor remainder at order $\omega_+(g_1)+1$ based at the point $\widetilde{y}_1$. This leads to an estimate of the (uniform) size of $\big(\textrm{Id}-\mathscr{T}_{g_1}\big){\bf g}_2$ by the converging sum $\sum_{m\geq 0} 2^{-m\{\omega(g_1)+\omega_+(g_1)+1\}}$. Note that the order $\omega_+(g_1)$ of expansion of $({\bf g}_2\backslash {\bf g}_1)$ in $\mathscr{T}_{g_1}$ does not depend on the scale assignment $\mu$.
\end{itemize}

{\it So, where does the {\sf BPHZ} formula \eqref{EqBPHZ} come from?} Given a scale assignment $\mu$ denote by ${\sf I}^{\pm}(\mu)$ the maximum and minimum of the scale assignments, respectively. For ${\sf I}^-(\mu)\leq i\leq {\sf I}^+(\mu)$ denote by $\bigsqcup_{\ell_i=1}^{n_i} g_{i,\ell_i}$ the connected components of the subset of $G$ which collects all the edges with scales bigger than or equal to $i$. One has for instance $n_{{\sf I}^-(\mu)}=1$ and $g_{n_{{\sf I}^-(\mu)},1} = G$. Two different $g_{i,\ell_i}, g_{j,\ell_j}$ are either disjoint or one is contained in the other. So this set of graphs is naturally organised into a tree structure when $i$ varies; one talks of {\it Gallavotti-Nicol\'o tree}. If $g_{i,\ell_i} \subset g_{j,\ell_j}$ then $g_{i,\ell_i}$ is quasi-local for $g_{j,\ell_j}$ for the assignment $\mu$. We define the parcimonious renormalization of ${\bf G}^\mu$ as
\begin{equation}
\label{EqFormulaParcimoniousRenormalisation}
^p\overline{\bf G}^\mu \defeq \prod_{i= {\sf I}^+(\mu)}^{{\sf I}^-(\mu)} \prod_{\ell_i=1}^{n_i} \big( \textrm{Id} - \mathscr{T}_{g_{i,\ell_i}} \big) {\bf G}^\mu.
\end{equation}
Because $\mathscr{T}_{g_{i,\ell_i}}{\bf G}^\mu=0$ if $g_{i,\ell_i}$ is not divergent, the parcimonious renormalization actually only involves the divergent sub-graphs of $G$. We order the product over $i$ from ${\sf I}^+(\mu)$ to ${\sf I}^-(\mu)$ to make it plaint that we do the product backward starting with the largest value of $i$. With the convention of Section \ref{SubsubsectionBPHZ} for the notation $\prod_{g\in\mathcal{F}}$ one can write
\begin{equation} \label{EqRenormAmu}
^p\overline{\bf G}^\mu = \sum_{\mathcal{F}\in\mathscr{F}(\mu)} \prod_{g\in \mathcal{F}} (-\mathscr{T}_g){\bf G}^\mu,
\end{equation}
for some particular forest $\mathscr{F}(\mu)$ of divergent sub-graphs of $G$. A refinement of the analysis of point (b) above is described in Section \ref{SubsubsectionForests} and shows that \textit{the summation of the $^p\overline{\bf G}^\mu$ over all the scale assignments $\mu$ is finite} -- Section \ref{SubsubsectionForests} gives some insights on that point. This is the main selling argument of the multiscale analysis: for each assignment the renormalization procedure specified by \eqref{EqFormulaParcimoniousRenormalisation} is actually well ordered and inductive. No renormalization operation made earlier is perturbed by a later operation -- so there is no problem with the so-called overlapping divergences in this picture. We define the \textbf{parcimonious renormalization of} $\bf G$ as 
\begin{equation} \label{EqParcimoniousRenormalisation}
\begin{split}
^p\overline{\bf G} &\defeq \sum_{\mu} {}^p\overline{\bf G}^\mu = \sum_{\mathcal{F}\in\mathscr{F}}\; \sum_{\mu; \mathcal{F}\in\mathscr{F}(\mu)} \prod_{g\in \mathcal{F}} (-\mathscr{T}_g){\bf G} ^\mu.
\end{split}
\end{equation}
Only the forests of quasi-local sub-diagrams are used here and renormalised via \eqref{EqFormulaParcimoniousRenormalisation}; they depend on each scale assignment $\mu$. We remark in particular that no overlapping divergent graphs are considered here. The formula defining $^p\overline{\bf G}$ has exactly the same form as the formula \eqref{EqBPHZ} definining $\overline{\bf G}$ except that in \eqref{EqBPHZ}  we sum over the set of all scale assignments, independently of $\mathcal{F}\in\mathscr{F}$ while we sum over the $\mu$ with $\mathcal{F}\in\mathscr{F}(\mu)$ in \eqref{EqParcimoniousRenormalisation}. One then sees that the BPHZ formula for $\overline{\bf G}$ comprises both the useful renormalizations from the parcimonious renormalization \eqref{EqFormulaParcimoniousRenormalisation} but also a number of operations that seem useless from the point of view of renormalization, at any fixed scale assignment $\mu$.

%%%-------------------------------------%%%
\subsubsection{Classification of forests$\boldmath{.}$}
\label{SubsubsectionForests}
%%%-------------------------------------%%%

Pick a forest of divergent sub-graphs $\mathcal{F}\in\mathscr{F}$. A \textit{sub-graph} $g$ of $G$, divergent or not, is said to be \textit{compatible with $\mathcal{F}$} if $\mathcal{F}\cup\{g\}$ is also a forest. One denotes then by $g^-_\mathcal{F}$ the unique element of $\mathcal{F}\cup\{g\}$ that strictly contains $g$ if $g$ is not the root of the forest, otherwise we set $g_\mcF^-=G$. We also denote by $g^+_\mathcal{F}$ the union of the elements of $\mathcal{F}$ that are strictly included in $g$. Given a scale assignment $\mu$, set then
\begin{equation*}
\begin{split}
m_{\mathcal{F},\mu}(g) &\defeq \min\Big\{\mu(e)\,;\,e\in \mathcal{E}\big(g\backslash g^+_\mathcal{F}\big)\Big\},   \\
M_{\mathcal{F},\mu}(g) &\defeq \max\Big\{\mu(e)\,;\, e_+\in \mathcal{V}(g)\,,\, e_-\in \mathcal{V}\big(g^-_\mathcal{F}\backslash g\big)\Big\}.
\end{split}
\end{equation*}
Define the {\it dangerous part of $\mathcal{F}$ for $\mu$} as the set
$$
\mathcal{D}_\mu(\mathcal{F}) \defeq \Big\{g\in \mathcal{F}\,;\, m_{\mathcal{F},\mu}(g) > M_{\mathcal{F},\mu}(g)\Big\}.
$$
The complement 
$$
\mathcal{S}_\mu(\mathcal{F})\defeq\mathcal{F}\backslash\mathcal{D}_\mu(\mathcal{F})
$$ 
of the dangerous part of $\mathcal{F}$ turns out to be a forest, called the {\it $\mu$-safe part of the forest $\mathcal{F}$}. A forest $\mcF$ for which $\mcF=\mcS_\mu(\mcF)$ is called a $\mu$-safe forest. One shows that
$$
\mathcal{S}_\mu\big(\mathcal{S}_\mu(\mathcal{F})\big) = \mathcal{S}_\mu(\mathcal{F}).
$$
This fact allows to partition $\mathscr{F}$ according to the $\mu$-safe part of its elements $\mcF$. We associate to $\mathcal{F}\in\mathscr{F}$ the following family of divergent sub-graphs of $G$ 
$$
\mathcal{S}^+_\mu(\mathcal{F}) \defeq \Big\{ g\subset G \,;\, g \textrm{ is compatible with }\mathcal{F}, \textrm{ divergent, and } g\in\mathcal{D}_\mu\big(\mathcal{S}_\mu(\mathcal{F})\cup\{g\}\big)\Big\}.
$$
The next statement was first proved by Feldman, Magnen, Rivasseau \& S\'en\'eor in \cite{FMRS} -- Lemma 2.3 therein; it is called the lemma of {\it classification of forests}.

\medskip

\begin{lem} \label{LemClassificationForests}
Let $\mcF\in\mathscr{F}$ be $\mu$-safe forest. Then $\mathcal{F}\cup \mathcal{S}^+_\mu(\mathcal{F}) = \mathcal{S}_\mu(\mathcal{F})\cup \mathcal{S}^+_\mu(\mathcal{F})$ is a forest. Another forest $\mathcal{F}'\in\mathscr{F}$ have the same $\mu$-safe part as $\mcF$ if and only if
$$
\mathcal{S}_\mu(\mathcal{F}) \subset \mathcal{F}' \subset \Big\{\mathcal{S}_\mu(\mathcal{F})\cup \mathcal{S}^+_\mu(\mathcal{F})\Big\}.
$$
\end{lem}

The set $\mathcal{S}_\mu(\mathcal{F})\cup \mathcal{S}^+_\mu(\mathcal{F})$ is in particular the maximal forest with the same $\mu$-safe part as the $\mu$-safe forest $\mathcal{F}$. Lemma \ref{LemClassificationForests} allows to rewrite $\overline{\bf G}$ as
\begin{equation}\label{Proof:EqPartitionForests}
\begin{split}
\overline{\bf G} &= \sum_{\mathcal{F}\in\mathscr{F}} \prod_{g\in \mathcal{F}} (-\mathscr{T}_g) {\bf G}   \\
&= \sum_\mu \sum_{\mathcal{F}\in\mathscr{F}} \prod_{g\in \mathcal{F}} (-\mathscr{T}_g){\bf G}^\mu   \\
&= \sum_\mu \;\sum_{\mathcal{G} \textrm{ safe for } \mu} \; \sum_{\mathcal{F};\mathcal{S}_\mu(\mathcal{F})=\mathcal{G}} \, \prod_{g\in \mathcal{F}} (-\mathscr{T}_g){\bf G}^\mu   \\
&= \sum_\mathcal{G\in\mathscr{F}} \; \sum_{\mu; \mathcal{G}=\mathcal{S}_\mu(\mathcal{G})} \; \sum_{\mathcal{F};\mathcal{S}_\mu(\mathcal{F})=\mathcal{G}} \, \prod_{g\in \mathcal{F}} (-\mathscr{T}_g){\bf G}^\mu   \\
&\hspace{-0.42cm}\overset{\textrm{Lemma \ref{LemClassificationForests}}}{=} \sum_\mathcal{G\in\mathscr{F}} \; \sum_{\mu; \mathcal{G}=\mathcal{S}_\mu(\mathcal{G})} \; \sum_{\mathcal{H}\subset \mathcal{S}^+_\mu(\mathcal{G})} \, \prod_{g\in \mathcal{G}\cup\mathcal{H}} (-\mathscr{T}_g){\bf G}^\mu.
\end{split}
\end{equation}
The term $\mathcal{G}=\emptyset$ in the last expression is equal to $^p\overline{\bf G}$. The set $\mathcal{S}_\mu(\mathcal{G})\cup\mathcal{S}^+_\mu(\mathcal{G})$ being a forest, it can be described as a union of disjoint trees $\mathcal{T}_1(\mathcal{G},\mu)$, $\dots$, $\mathcal{T}_{k(\mathcal{G}, \mu)}(\mathcal{G},\mu)$. The index $k(\mathcal{G},\mu)$ is bounded above by a constant independent of $\mathcal{G}$ and $\mu$. One can then write
\begin{equation} \label{EqPartBPHZFormula}
\sum_{\mathcal{H}\subset \mathcal{S}^+_\mu(\mathcal{G})} \prod_{g\in \mathcal{G}\cup\mathcal{H}} (-\mathscr{T}_g){\bf G}^\mu = \prod_{k=1}^{k(\mathcal{G},\mu)} \prod_{g\in \mathcal{T}_k(\mathcal{G},\mu)} \mathscr{R}_g {\bf G}^\mu,
\end{equation}
where
$$
\begin{array}{ll}
\mathscr{R}_g = -\mathscr{T}_g, &\textrm{if } g\in\mathcal{G},   \\
\mathscr{R}_g = \textrm{Id}-\mathscr{T}_g, & \textrm{if } g\in\mathcal{S}^+_\mu(\mathcal{G}).   \\
\end{array}
$$
The operators involved in different $\mathcal{T}_k(\mathcal{G},\mu)$ commute. The index set $\mathscr{F}$ in the fifth equality of \eqref{Proof:EqPartitionForests} for $\overline{\bf G}$ being finite, it suffices to consider the convergence of each sum $\sum_{\mu; \mathcal{S}_\mu(\mathcal{G})=\mathcal{G}}$. Let us concentrate on each product $\prod_{g\in \mathcal{T}_k(\mathcal{G},\mu)}$. Setting $-\mathscr{T}_G \defeq \textrm{Id}$, one can re-index the product on the tree $\mathcal{T}_k(\mathcal{G},\mu)$ in \eqref{EqPartBPHZFormula} and write
$$
\prod_{g\in \mathcal{T}_k(\mathcal{G},\mu)} \mathscr{R}_g = \prod_{g\in (\mcG\cap\mathcal{T}_k(\mathcal{G},\mu))\cup\{G\}} (-\mathscr{T}_g) \prod_{h\in(\mathcal{S}^+_\mu(\mathcal{G}) \cap\mathcal{T}_k(\mathcal{G},\mu)); h_\mathcal{G}^-=g} \big(\textrm{Id}-\mathscr{T}_h\big).
$$
Note that the two operators $\mathscr{T}_g$ and $\mathscr{T}_h$ commute in this formula since they act on the outside and the inside of $g$, respectively. The bounds on $\prod_{g\in \mathcal{T}_k(\mathcal{G},\mu)} \mathscr{R}_g {\bf G}^\mu$ that one can get depend on the order in which one does the integration. To make that point plain, look at the example 
$$
\int \prod_{i=1}^3 a_i^{-\frac{d}{2}} e^{-\frac{\vert x-y\vert^2}{a_i}} \prod_{j=1}^2 b_j^{-\frac{d}{2}} e^{-\frac{\vert y-z\vert^2}{b_j}} c^{-\frac{d}{2}} e^{-\frac{\vert z-x\vert^2}{c}}\,dxdydz
$$
that corresponds to some graph ${\sf G}= (\sf E,V)$ that you are invited to draw. One can choose to bound some of the exponentials by $1$ and only keep one term for each integration variable. This corresponds to choosing a covering tree $\sf T$ of $\sf G$. Denoting by $r_e\in\{a_1,a_2,a_3,b_1,b_2,c\}$ the parameter associated to the edge $e\in{\sf E}$ the choice of a covering tree $\sf T$ gives for the integral the following $\sf T$-dependent bound
$$
\left(\prod_{e\in {\sf E}} r_e\right)^{-\frac{d}{2}} \prod_{e\in{\sf T}} r_e^{\frac{d}{2}}
$$
that one can optimize by an appropriate choice of $\sf T$ guided by the form of the bound. Back to our general setting, recall that we denote for each $i$ by $(g/g^+_\mathcal{G})_{(i,\ell_i)}$ the set of connected component of the graph $g/g^+_\mathcal{G}$ with scale bigger than or equal to $i$. We choose a covering tree of $G$ whose restriction to each $(g/g^+_\mathcal{G})_{(i,\ell_i)}$ is a sub-tree. For these graphs, the equivalent of the previous $r_e$ are essentially given by $2^{-i}$. This way, we get as many  factors $2^{-id/2}$ as we can get from the integration variables in the $(g/g^+_\mathcal{G})_{(i,\ell_i)}$ if the kernels that we keep are untouched. If a connected component $(g/g^+_\mathcal{G})_{(i,\ell_i)}$ is not quasi-local, and if $(g/g^+_\mathcal{G})_{(i,\ell_i)} \neq g/g^+_\mathcal{G}$, we have $\omega((g/g^+_\mathcal{G})_{(i,\ell_i)})\geq 1$ and it gives a factor $2^{-\omega((g/g^+_\mathcal{G})_{(i,\ell_i)})}$; if $(g/g^+_\mathcal{G})_{(i,\ell_i)} = g/g^+_\mathcal{G}$, one has a factor $O(1)$. If the connected component $(g/g^+_\mathcal{G})_{(i,\ell_i)}$ is divergent and quasi-local for $\mu$ it is part of the collection of  $\big\{ h\in\mathcal{S}^+_\mu(\mathcal{G}); h_\mathcal{G}^-=g \big\}$, and the operator $\big(\textrm{Id}-\mathscr{T}_h\big)$ gives an ad hoc factor that provides a total contribution $O(1)$. To summarize, set
$$
\omega'\big((g/g^+_\mathcal{G})_{(i,\ell_i)}\big) := \Big(1\vee \omega\big((g/g^+_\mathcal{G})_{(i,\ell_i)}\big)\Big)\,{\bf 1}_{g/g^+_\mathcal{G} \neq (g/g^+_\mathcal{G})_{(i,\ell_i)}}.
$$ 
The above analysis leads to an estimate of the form
$$
\bigg| \prod_{g\in \mathcal{T}_k(\mathcal{G},\mu)} \mathscr{R}_g {\bf G}^\mu\bigg| \lesssim \prod_{g\in \mathcal{T}_k(\mathcal{G},\mu)\cup\{G\}} \, \prod_{(i,\ell_i)} 2^{-\omega'((g/g^+_\mathcal{G})_{(i,\ell_i)})}.
$$
Note that the scale assignment implicitly appears in the right hand side in the bounds for $i$. It is then possible from this bound to obtain the summability in $\mu$ of \eqref{EqPartBPHZFormula}, hence of $\overline{\bf G}$. See Section 3.2 of \cite{HairerBPHZ} for the model case of Feynman graphs and Section 8 of \cite{CH} for the regularity structure case. (Our account of the BPHZ procedure was influenced by the very nice presentation given by F. Vignes-Tourneret in his PhD thesis \cite{VignesTourneret}.)

\medskip

%%%-----------------------%%%
\subsubsection{Back to trees$\boldmath{.}$}
\label{SubsubsectionBackToTrees}
%%%-----------------------%%%

The preceding sections sketch the main points of the mechanics involved in the BPHZ renormalization of some iterated integrals as $\bf G$ or \eqref{EqL2OmegaNormModel}. However, in a regularity structures setting, we do not directly have a hand on \eqref{EqL2OmegaNormModel} but rather on the $\overline{\bf W}^{\epsilon,m}_\tau$, via the definition of the renormalized $\overline{\sf \Pi}_0^\epsilon\tau$. The renormalization procedure introduced by Bruned, Hairer \& Zambotti in \cite{BHZ} has this property that it somehow induces on the quantities \eqref{EqL2OmegaNormModel} built from the renormalized model the same good renormalization/(Taylor remainder) operations as in the BPHZ procedure. This is made possible by the fact that in the $m$-th term of the pairing \eqref{EqL2OmegaNormModel}, except from the whole graph itself, every divergent sub-graph is either in $\overline{\bf W}^{(\epsilon_1,\epsilon_2),m}_\tau(a,y)$ or in $\overline{\bf W}^{(\epsilon_1,\epsilon_2),m}_\tau(b,y)$: it does not overlap both parts at a time. The appropriate definition of the renormalized smooth model $\overline{\sf \Pi}^\epsilon$ was first given in Section 6 of \cite{BHZ}.

\medskip

%%-----------------------%%
\subsection{Recentering$\boldmath{.}$ \hspace{0.1cm}}
\label{SubsectionRecentering}
%%-----------------------%%

Producing some $(\epsilon_1,\epsilon_2)$-independent finite bounds for \eqref{EqL2OmegaNormModel} is not sufficient to prove that the renormalized smooth models $\overline{\sf \Pi}^\epsilon$ converge to a limit in the space of models. We need for that purpose the bounds to be of the form $o_{\max(\epsilon_1,\epsilon_2)}(1) \, \lambda^{2(r(\tau)+\kappa)}$ as in \eqref{EqEstimateConvergenceModels}. The BPHZ mechanics explained above does not help in that task and provides no clue on that point for the reason that it essentially deals with the non-recentered interpretation map $\overline{\bf \Pi}^\varepsilon$ rather than with the $\overline{\Pi}_x^\varepsilon$. A different type of argument is needed to understand this scaling behaviour obtained from recentering. Theorem A.3 in Hairer \& Quastel's work \cite{HQ}  is a good entry point for that question. It was latter improved by Bruned \& Nadeem in Theorem 3.1 of \cite{BrunedNadeem}. The latter was used in \cite{BailleulBrunedGKPZ} by Bailleul \& Bruned to give a simple proof of convergence of the BHZ renormalized model of the one dimensional generalized (KPZ) equation driven by a spacetime white noise.

\medskip

%%-----------------------------------------------------------------------------------------------------------%
%\section{Linares, Otto, Tempelmayr \& Tsatsoulis and noise derivative}% on multi-indices}
%\label{SectionLOTT}
%%-----------------------------------------------------------------------------------------------------------%

%\vfill \pagebreak

%------------------------------------------%
\section{Non-diagramatic methods}
\label{SectionNondiagramatic}
%------------------------------------------%

Chandra \& Hairer assume in \cite{CH} some moment type condition on the law of the noise in the form of some quantitative estimates on its cumulant. Within that setting the convergence criterion \eqref{EqEstimateConvergenceModels} is relatively sharp and one does not loose much in quantifying the expected convergence result of Theorem \ref{HS:thm:main} in the form \eqref{EqEstimateConvergenceModels}. The strategy of \cite{CH} is somewhat optimal: no information is lost in the chaos expansion and each iterated integral that comes with that decomposition is proved to scale as $o_{\max(\epsilon_1,\epsilon_2)}(1)\,\lambda^{2(r(\tau)+\kappa)}$. Linares, Otto, Tempelmayr \& Tsatsoulis introduced in \cite{LOTT} a different type of strategy in their development of an alternative to regularity structures well suited for the study of a class of quasilinear singular SPDEs. A probability measure $\bbP$ satisfies a spectral gap inequality if the variance of any $L^2(\bbP)$ random variable is bounded above by a constant multiple of the variance of its Malliavin derivative. Trading Chandra \& Hairer's moment assumption for a spectral gap assumption on the law of the noise opens the door to a different kind of strategy for proving the convergence result of Theorem \ref{HS:thm:main}. Work in a regularity structure with a noise symbol and a noise derivative symbol, with some models where one can represent the $\Pi_x\tau$ for the usual symbols $\tau$ but also their Malliavin derivative. The reasoning is inductive and informally expressed as follows.
\begin{itemize}\setlength{\itemsep}{0.1cm}
	\item[--] Use the spectral gap inequality to propagate a convergence result for the noise derivative of $\Pi_x\tau$ to $\Pi_x\tau$ itself.
	
	\item[--] Show that the convergence of a number of $\Pi_x\tau$ and their derivatives allow to prove the convergence of the noise derivative of $\Pi_x\sigma$ for a new symbol $\sigma$ not in the list of the previous symbols $\tau$.
\end{itemize}
We choose to leave aside in this review the very important work \cite{LOTT} of Linares, Otto, Tempelmayr \& Tsatsoulis and concentrate on the works \cite{HS, RandomModel} that implement the strategy of \cite{LOTT} in a regularity structure setting indexed by trees, as they share the same technical background. This makes their introduction and comparison easier. For a reader interested in the multi-indices approach to regularity structures we recommend the lecture notes \cite{LO, BOT}. The works \cite{HS} and \cite{RandomModel} do not implement the above inductive strategy in the same way. We describe the details of \cite{HS} and \cite{RandomModel} in Section \ref{SectionHS} and Section \ref{SectionBH}, respectively. The contents of these sections are independent from each other; they only share some basic matters related to the Malliavin derivative  (Section \ref{SubsectionHS1}) and the initial setup for the proof described at the beginning of Section \ref{SubsectionHS3}.

\medskip

%--------------------------------------------------------------------%
\subsection{Hairer \& Steele \cite{HS} and pointed modelled distributions$\boldmath{.}$ \hspace{0.1cm}}
\label{SectionHS}
%--------------------------------------------------------------------%

We stated in Theorem \ref{HS:thm:main} the main convergence result of \cite{HS} for the BPHZ renormalized models. We describe in this section the main ideas involved in the proof of this statement.

\bigskip

%%-----------------------------------------------------------------%%
\subsubsection{Poincar\'e inequality and Malliavin derivative$\boldmath{.}$ \hspace{0.1cm}}
\label{SubsectionHS1}
%%-----------------------------------------------------------------%%

Following the insight of \cite{LOTT} the key tool used in \cite{HS} is the {\bf Poincar\'e inequality}. It is satisfied by the law of the white noise. Denote by $H$ the Hilbert space $L^2(\bbR^d)$. A cylindrical function $F:\mcD'(\bbR^d)\to\bbR$ is a function of the form
$$
F(\xi) = f\big(\xi(\varphi_1),\dots,\xi(\varphi_N)\big),
$$
where $f\in C^\infty(\bbR^N)$ has at most polynomial growth and $\varphi_1,\dots,\varphi_N\in\mcD(\bbR^d)$. For each cylindrical function $F$ the Malliavin derivative in the direction $h\in H$ is defined by
$$
\nabla_hF(\xi)\defeq
\lim_{\varepsilon\to0}\frac{F(\xi+\varepsilon h)-F(\xi)}\varepsilon=
\sum_{i=1}^N\partial_if\big(\xi(\varphi_1),\dots,\xi(\varphi_N)\big) \, \langle h,\varphi_i\rangle_{L^2}.
$$
We denote by $\nabla F(\xi)\in H^*$ the linear form $h\mapsto\nabla_hF(\xi)$. Denoting by $\bbP$ the probability law of the white noise and by $\bbE$ the expectation with respect to $\bbP$, we define $\bbD^{1,2}$ as the completion of the set of all cylindrical functions under the norm $\big(\bbE[\|F\|^2]+\bbE[\|\nabla F\|_{H^*}^2]\big)^{1/2}$.
It is well known that $\bbP$ satisfies the Poincar\'e inequality
\begin{equation} \label{HS:EqSGInequality}
\bbE\Big[ \big( F-\bbE[F] \big)^2\Big] \leq \bbE\big[\|\nabla F\|_{H^*}^2\big]
\end{equation}
for any $F\in\bbD^{1,2}$. One can find the classical proof based on the log-Sobolev inequality in Theorem 5.5.1 of Bogachev's book \cite{Bogachev}. The inequality \eqref{HS:EqSGInequality} is also called a {\it spectral gap inequality} as it turns out to be a consequence of an estimate for the spectral gap of the Ornstein-Uhlenbeck operator defined on $L^2(\Omega)$. The inequality \eqref{HS:EqSGInequality} allows to reduce the task of estimating the $L^2$ size of a polynomial functional of  a white noise to the $L^2$ size of another functional with lower degree. In our setting we can write for any $\tau\in\mcB$ and any $q\ge1$
\begin{align} \label{HS:useofpoincare}
\bbE\big[|(\overline{\Pi}^\epsilon_x\tau)(\varphi_x^\lambda)|^q\big] 
&\lesssim \big|\bbE\big[(\overline{\Pi}^\epsilon_x\tau)(\varphi_x^\lambda)\big]\big|^q 
+ \bbE\bigg[\sup_{\|h\|_{H}=1}|\nabla_h(\overline{\Pi}_x^\epsilon\tau)(\varphi_x^\lambda)|^q\bigg].
\end{align}
As stated in the end of Section \ref{SubsectionBasicsBPHZ}, the BPHZ renormalized model $\overline{\sf M}^\epsilon$ is the unique model associated with some particular type of preparation maps with the property that $\bbE[(\overline{\bf \Pi}^\epsilon\tau)(x)]=0$ everywhere, for all $\tau\in\mcB$ with negative degree. This property provides at some relatively low cost some estimates on $\bbE\big[(\overline{\Pi}^\epsilon_x\tau)(\varphi_x^\lambda)\big]$ of the right order in terms of $\lambda$ --  see the proof of Proposition 5.2 of \cite{HS}. The translation invariance of both the kernel and the law of the noise are essentially only used only at that point. (Precisely, Hairer and Steele first proved the convergence of a modified version $\widetilde{\sf M}^\varepsilon$ defined by the condition $\bbE[(\widetilde{\Pi}_0^\epsilon\tau)(\varphi)]=0$ for some fixed test function $\varphi$, and called $\overline{\text{BPHZ}}$ model therein in \cite{HS}, and then deduced the convergence of $\overline{\sf M}^\varepsilon$ from the convergence of $\widetilde{\sf M}^\varepsilon$.) The spectral gap inequality \eqref{HS:useofpoincare} then brings back the analysis to estimating the Malliavin derivative term. We indicated above that this is the key fact for an inductive approach to the convergence of ${\overline{\sf M}^\epsilon}$.

\ssk

It seems natural to introduce a new symbol $\dot{\Xi}$ representing a generic element of $H$ and to interpret the derivative $\nabla_h(\overline{\Pi}_x^\epsilon\tau)$ as the application of (some extension of the interpretation map) $\overline{\Pi}_x^\epsilon$ to some symbol involving $\dot{\Xi}$. Considering each symbol $\tau\in\mcF$ as a multilinear functional of $\Xi$, we define $\dot{\mcF}$ as the set of symbols obtained by replacing by $\dot{\Xi}$ any of the arguments $\Xi$ in these functions. For example, from the symbol $\Xi\mcI(\Xi)^2$ in $\mcF$ we obtain the two symbols $\dot{\Xi}\mcI(\Xi)^2$ and $\Xi\mcI(\Xi)\mcI(\dot{\Xi})$ in $\dot{\mcF}$. By using these additional symbols we can define the formal Malliavin derivative operator $D:\spa(\mcF)\to\spa(\dot{\mcF})$ as follows:
\begin{align*}
D\Xi=\dot{\Xi},\qquad
DX^k=0,\qquad
D\mcI_k(\tau)=\mcI_k(D\tau),\qquad
D(\tau\sigma)=(D\tau)\sigma+\tau(D\sigma).
\end{align*}
(The map $\mcI_k$ is extended linearly in the third equality). Let now $\dot{\mcB}$ be the subset of $\dot{\mcF}$ consisting of all the symbols obtained by replacing any one argument $\Xi$ in $\tau\in\mcB$ with $\dot{\Xi}$. Set
$$
\widetilde{\bfT}\defeq\spa(\mcB\cup\dot{\mcB}).
$$ 
In \cite{HS} the authors extended the degree map to $\widetilde{\bfT}$ setting 
$$
r(\dot{\Xi})=\theta,\quad \textrm{ for some fixed } \theta\in(0,1).
$$ 
They also constructed a regularity structure $\widetilde{\scT}$ with model space $\widetilde{\bfT}$. Using an argument based on the reconstruction theorem (see also Lemma \ref{HS:lem:reconst} below), they were able to extend the BPHZ model $\overline{\sf M}^\epsilon$ into a unique model $\widetilde{\sf M}^{\epsilon,h}$ on $\widetilde{\scT}$ satisfying $\widetilde{\bf\Pi}_x^{\epsilon,h}(\dot{\Xi}) = h^\epsilon = h*\varrho^\varepsilon$. One might expect that the identity
\begin{equation} \label{EqIntertwiningNablaD}
\nabla_h\widetilde{\Pi}_x^{\epsilon}(\tau) = \widetilde{\Pi}_x^{\epsilon,h}(D\tau)
\end{equation}
holds for any $\tau\in\mcB$, but this is not the case, even for $\tau=\Xi$ for which we have
$$
\nabla_h\overline{\Pi}_x^{\epsilon}(\Xi) = \nabla_h\xi^\epsilon = h^\epsilon,\qquad
\widetilde{\Pi}_x^{\epsilon,h}(D\Xi) = \widetilde{\Pi}_x^{\epsilon,h}(\dot{\Xi}) = h^\epsilon-h^\epsilon(x).
$$
This difference arises because $\dot{\Xi}$ has a much higher degree than $\Xi$. To circumvent this problem Hairer \& Steele defined some $x$-dependent modelled distributions $\widetilde{f}_x^{\epsilon;\tau}$ satisfying the intertwining relation
\begin{align} \label{HS:eq:realandformalderivative}
\nabla_h\overline{\Pi}_x^{\epsilon}(\tau) = \widetilde{\Pi}_x^{\epsilon,h}\big(\widetilde{f}_x^{\epsilon;\tau}(x)\big) = \widetilde\mcR^{\epsilon,h}\big(\widetilde{f}_x^{\epsilon;\tau}\big).
\end{align}
To capture the main analytic properties of these functions $\widetilde{f}_x^{\varepsilon;\tau}$ they introduced in Section 3 of \cite{HS} a new notion of modelled distribution that is the object of the next section. We return to the intertwining relation in Section \ref{SectionBH}.

\bigskip

%%----------------------------------------------%%
\subsubsection{Pointed modelled distributions$\boldmath{.}$ \hspace{0.1cm}}
\label{SubsectionHS2}
%%----------------------------------------------%%

In this subsection we let $\scT=(\bfA,\bfT,\bfG)$ be an arbitrary regularity structure with regularity $\alpha_0$ and let ${\sf M} = (\Pi,\Gamma)$ be a model for $\scT$ on $\bbR^d$. First, we recall the definition of some $B_{p,\infty}^\gamma$-type space of modelled distributions. For any $\gamma\in\bbR$ and $p\in[1,\infty]$ we define $\mcD_p^\gamma=\mcD_p^\gamma(\Gamma)$ as the space of all functions $f:\bbR^d\to\bfT_{<\gamma}$ for which
\begin{align*}
\rest{f}_{p,\gamma;C}&\defeq
\max_{\alpha\in\bfA,\,\alpha<\gamma}\big\|\|f(x)\|_{\alpha}\big\|_{L_x^p(C)}<\infty,
\\
\|f\|_{p,\gamma;C}&\defeq
\max_{\alpha\in\bfA,\,\alpha<\gamma}
\sup_{0<\|y\|\le1}
\frac{\big\| \big\| f(x+y) - \Gamma_{(x+y)x}(f(x)) \big\|_\alpha\big\|_{L_x^p(C)}}{\|y\|^{\gamma-\alpha}}<\infty
\end{align*}
for any compact subset $C\subset\bbR^d$. We define the notations $\tri f\tri_{p,\gamma;C}$ and $\tri f^{(1)},f^{(2)}\tri_{p,\gamma;C}$ in a way similar to Definition \ref{HS:def:MD}. There is a reconstruction theorem for this kind of modelled distributions which reads as follows: For any $\gamma>0$ and $p\in[1,\infty]$ there exists a unique continuous linear map $\mcR^{\sf M}:\mcD_p^\gamma(\Gamma)\to\mcD'(\bbR^d)$ satisfying
$$
\sup_{\varphi\in\mcB_r}\sup_{\lambda\in(0,1]}\lambda^{-\gamma}
\big\|\big(\mcR^{\sf M} f-\Pi_xf(x)\big)(\varphi_x^\lambda)\big\|_{L_x^p(C)}
\lesssim\|\Pi\|_{\gamma;\overline{C}}\|f\|_{p,\gamma;\overline{C}}.
$$
Moreover the mapping $({\sf M},f)\mapsto\mcR^{\sf M} f$ is locally Lipschitz continuous with respect to the quasi-metrics $\tri {\sf M}^{(1)} \,,\, {\sf M}^{(2)}\tri_{\gamma;\overline{C}}$ and $\tri f^{(1)},f^{(2)}\tri_{\gamma;\overline{C}}$. (This is the statement of Theorem 3.4 in \cite{HS}.)

\medskip

Given any $x\in\bbR^d$, an $x$-{\bf pointed modelled distribution} is a modelled distribution in the class $\mcD_p^\gamma$ which behaves better near $x$ than it does elsewhere. Here is the archetypal example. For any smooth function $f:\bbR^d\to\bbR$ and any integers $\gamma<\nu$ we define the function $F_x:\bbR^d\to\spa\{X^k\}_{k\in\bbN^d}$ by
$$
F_x(y) \defeq \sum_{|k|<\gamma}\bigg(\partial^kf(y)-\sum_{|\ell|<\nu-|k|}\frac{\partial^{k+\ell}f(x)}{\ell!}(y-x)^\ell\bigg)\frac{X^k}{k!} \eqdef \sum_{|k|<\gamma}\partial^k f_x(y) \, \frac{X^k}{k!}.
$$
We read from Taylor theorem the identity
\begin{align*}
F_x(z)-\Gamma_{zy}F_x(y)
=\sum_{|k|<\gamma}\frac{X^k}{k!}\sum_{|\ell|=\gamma-|k|}\frac{(z-y)^\ell}{\ell!}
\int_0^1|\ell|(1-t)^{|\ell|}\partial^{k+\ell}f_x\big(y+t(z-y)\big)dt,
\end{align*}
so we have
$$
\big\| F_x(z)-\Gamma_{zy}F_x(y) \big\|_{|k|} \lesssim \|z-y\|^{\gamma-|k|}
$$
globally in $y,z$. However since $|\partial^{k+\ell}f_x(y)|\lesssim\|y-x\|^{\nu-|k+\ell|}$ the above estimate improves to
$$
\big\| F_x(z)-\Gamma_{zy}(F_x(y)) \big\|_{|k|} \lesssim \lambda^{\nu-\gamma}\|z-y\|^{\gamma-|k|}
$$
for $y$ and $z$ in the ball $B(x,\lambda)$ of center $x$ and radius $\lambda>0$. For any $\gamma,\nu\in\bbR$, $p\in[1,\infty]$, and $x\in\bbR^d$ we define
$$
\mcD_p^{\gamma,\nu;x}=\mcD_p^{\gamma,\nu;x}(\Gamma)
$$ 
as the space of all functions $f\in\mcD_p^\gamma(\Gamma)$ for which
\begin{align*}
\rest{f}_{p,\gamma,\nu;x}&\defeq
\max_{\alpha\in\bfA,\,\alpha<\gamma}\sup_{\lambda\in(0,1]}
\lambda^{\alpha-\nu}\big\|\|f(y)\|_{\alpha}\big\|_{L_y^p(B(x,\lambda))}<\infty,
\\
\|f\|_{p,\gamma,\nu;x}&\defeq
\max_{\alpha\in\bfA,\,\alpha<\gamma}\sup_{\lambda\in(0,1]}
\sup_{0<\|z\|\le\lambda}
\lambda^{\gamma-\nu}
\frac{ \big\| \big\| f(y+z)-\Gamma_{(y+z)y}(f(y)) \big\|_\alpha \big\|_{L_y^p(B(x,\lambda))}}{\|z\|^{\gamma-\alpha}}<\infty.
\end{align*}
The reconstruction of a pointed modelled distribution satisfies an improved version of the reconstruction estimate that involves the following quantity. Fix a positive integer $r>|\alpha_0|$. For $f\in\mcD_p^{\gamma,\nu;x}(\Gamma)$ we set
$$
m\big(f,\Lambda;x\big)
\defeq
\sup_{\varphi\in\mcB_r}\bigg\{\sup_{\lambda\in(0,1]}\lambda^{d/p-\nu}|\Lambda(\varphi_x^\lambda)| + \sup_{0<\delta\le\lambda\le1}
\lambda^{\gamma-\nu}\delta^{-\gamma}\big\|\big(\Lambda-\Pi_yf(y)\big)(\varphi_y^\delta)\big\|_{L_y^p(B(x,\lambda))}\bigg\}
$$

\medskip

\begin{thm} \label{HS:thm:pointreconst}
\cite[Theorem 3.15]{HS}
Let $p\in[1,\infty]$, $\gamma\in(0,\alpha_0+d/p)\setminus\bbN$, $\nu\in\bbR$, and $x\in\bbR^d$. For any $f\in\mcD_p^{\gamma,\nu;x}(\Gamma)$ its reconstruction $\mcR^{\sf M} f$ satisfies the estimate
$$
m\big(f,\mcR^{\sf M} f;x\big) \lesssim
\tri f\tri_{p,\gamma,\nu;x}\|\Pi\|_{\gamma;B(x,2)} \big( 1+\|\Gamma\|_{\gamma;B(x,2)} \big).
$$
\end{thm}

\medskip

%%----------------------------------------------------------------------%%
\subsubsection{A sketch of the proof of Theorem \ref{HS:thm:main}$\boldmath{.}$ \hspace{0.1cm}}
\label{SubsectionHS3}
%%----------------------------------------------------------------------%%

We describe in this section the architecture of the argument used by Hairer \& Steele in \cite{HS} to prove Theorem \ref{HS:thm:main}. We focus here of the mechanics to get some $\epsilon$-uniform bounds for $\bbE\big[\tri\overline{\sf M}^\epsilon\tri_{\gamma;C}^q\big]$. Especially, we concentrate on the bounds of the $\overline{\Pi}^\varepsilon$ part since the bounds of the $\overline{\Gamma}^\varepsilon$ part automatically follow -- see Theorems 5.14 and 10.7 of \cite{Hairer}. To show the convergence of $\bbE\big[\tri\overline{\sf M}^{\epsilon_1} \,,\, \overline{\sf M}^{\epsilon_2}\tri_{\gamma;C}^q\big]$ to $0$ as $\epsilon_1,\epsilon_2$ go to $0$ one needs to introduce a new notion of modelled distributions measured in some negative Sobolev norms. We do not touch upon this point here and refer the reader to Section 3.3 of \cite{HS}.

\ssk

The proof of the $\epsilon$-uniform bound is done by induction with respect to the following pre-order. Denote by $n_\Xi(\tau)$ the number of symbols $\Xi$ contained in $\tau$ and define the preorder $\preceq$ by setting
\begin{equation} \label{EqDefnOrder}
\sigma \preceq \tau \quad \overset{\text{def}}{\Longleftrightarrow} \quad \big(n_\Xi(\sigma), |E_\sigma|,r(\sigma)\big) \leq \big(n_\Xi(\tau),|E_\tau|,r(\tau)\big)
\end{equation}
with the inequality $\leq$ in the right hand side standing here for the lexicographical order. We order the elements of $\mcB\setminus\{X^k\}_{k\in\bbN^d}$ as
$$
\mcB\setminus\{X^{k}\}_{{k}\in\bbN^d} = \big\{\tau_1\preceq\tau_2\preceq\cdots\big\}
$$
and set for each $i\ge1$
$$
\mcB_i \defeq \big\{ \tau_1,\tau_2,\dots,\tau_i \big\}.
$$
We need some notations to keep track of the sizes of $\overline\Pi^\epsilon$ and $\overline\Gamma^\epsilon$ on some vector spaces associated with $\mcB_i$. We introduce for that purpose for each $\tau\in\mcB$ the quantities
\begin{align*}
\| \overline\Pi^\epsilon \|_{\tau;C} &\defeq 
\sup_{\varphi\in\mcB_r}\sup_{\lambda\in(0,1]}\sup_{x\in C}
 \lambda^{-r(\tau)} \big| \big(\overline\Pi^\epsilon_x\tau\big) (\varphi_x^\lambda) \big|,
\\
\| \overline\Gamma^\epsilon \|_{\tau;C} &\defeq
\max_{\alpha\in\bfA,\,\alpha<r(\tau)}
\sup_{x,x+y\in C,\,y\neq0}
\frac{\big\| \overline\Gamma^\epsilon_{(x+y)x}(\tau) \big\|_\alpha}
{\|y\|^{r(\tau)-\alpha}},   
\\
\tri \overline{\sf M}^\epsilon\tri_{\tau;C} &\defeq \| \overline\Pi^\epsilon \|_{\tau;C} + \| \overline\Gamma^\epsilon \|_{\tau;C}.
\end{align*}
For any finite subset $\mcA\subset\mcB$ we also write
\begin{align*}
\| \overline\Pi^\epsilon \|_{\mcA;C} &\defeq \max_{\tau\in\mcA} \| \overline\Pi^\epsilon \|_{\tau;C},
\\
\| \overline\Gamma^\epsilon \|_{\mcA;C} &\defeq \max_{\tau\in\mcA} \| \overline\Gamma^\epsilon \|_{\tau;C},
\\
\tri \overline{\sf M}^\epsilon\tri_{\mcA;C} &\defeq \| \overline\Pi^\epsilon \|_{\mcA;C} + \| \overline\Gamma^\epsilon \|_{\mcA;C}.
\end{align*}

\ssk

Assume now that we have proved some $\epsilon$-uniform bounds for the model on $\mcB_{i-1}$ and let us derive the estimates on $\mcB_i$. The estimate on the initial set $\mcB_0\defeq\emptyset$ is free. We use some different tools in the induction step depending on the sign of $r(\tau_i)$. As a shorthand notation we denote by $\overline{\mcR}^\epsilon$ the reconstruction operator associated with $\overline{\sf M}^\epsilon$. We write $\widetilde{\mcR}^{\epsilon,h}$ for the reconstruction operator associated with $\widetilde{\sf M}^{\epsilon,h}$, for any $h\in H$.
\begin{enumerate}\setlength{\itemsep}{0.1cm}
	\item[$\rhd$] If $r(\tau_i)>0$ the estimate follows from the reconstruction theorem applied to a well-chosen modelled distribution with positive regularity -- see Proposition 3.31 of \cite{Hairer}.
The key fact is that the function $\gamma_x^{\epsilon;\tau_i}:\bbR^d\to\bfT$ defined by
$$
\gamma_x^{\epsilon;\tau_i}(y)\defeq \overline\Gamma^\epsilon_{yx}(\tau_i) - \tau_i
$$
is a modelled distribution in the class $\mcD^{r(\tau_i)}$ associated with $\overline{\sf M}^\epsilon$. Its reconstruction is given by
$$
\overline\mcR^{\epsilon} \big(\gamma_x^{\epsilon;\tau_i} \big)(y) = \overline\Pi^\epsilon_y\big(\gamma_x^{\epsilon;\tau_i}(y)\big)(y) = \big( \overline\Pi^\epsilon_x\tau_i - \overline\Pi^\epsilon_y\tau_i \big)(y) = \big( \overline\Pi^\epsilon_x\tau_i\big)(y).
$$

\ssk

\begin{lem} \label{HS:lem:reconst}
If $r(\tau_i)>0$ there exists a polynomial $P$ such that for any compact set $C\subset\bbR^d$ we have almost surely $\|\overline{\Pi}^{\epsilon}\|_{\tau_i;C} \le P\big(\|\overline{\Pi}^{\epsilon}\|_{\mcB_{i-1};\overline{C}}\big)$.
\end{lem}

\ssk

The reasoning involved in this item is purely deterministic.   \vspace{0.1cm}
	
	\item[$\rhd$] If $r(\tau_i)\le0$ we use the spectral gap inequality \eqref{HS:useofpoincare} and focus on the task of obtaining some $\epsilon$-uniform bound on the Malliavin derivative term. This is done in two steps that are detailed below.

\begin{itemize} \setlength{\itemsep}{0.1cm}
\item[--]
{\bf Algebraic step}:
One constructs some $x$-dependent modelled distributions $\widetilde{f}_x^{\epsilon;\tau}$ satisfying the intertwining relation \eqref{HS:eq:realandformalderivative} with the Malliavin derivative operator. 

\item[--]
{\bf Analytic step}:
One shows that a certain truncation of $\widetilde{f}_x^{\epsilon;\tau}$ is a pointed modelled distribution that satisfies some good estimate.
\end{itemize}

\ssk

\noindent
As in Lemma \ref{HS:lem:reconst}, these two steps provide a probabilistic control of $\|\overline{\Pi}^{\epsilon}\|_{\tau_i;C}$ in terms of $\|\overline{\Pi}^{\epsilon}\|_{\mcB_{i-1};\overline{C}}$, which is sufficient to close the induction in the case where $r(\tau_i)\le0$. The reasoning involved in this item is probabilistic as it rests on the spectral gap inequality.
\end{enumerate}

\medskip

%%-----------------------%%
\noindent -- \textbf{Algebraic step$\boldmath{.}$}
%%-----------------------%%
For any $\tau\in\mcB$ we define the function $f_x^{\epsilon;\tau} : \bbR^d\to\bfT$ by
$$
f_x^{\epsilon;\tau}(y) \defeq \overline\Gamma^\epsilon_{yx}(\tau).
$$
Since one has for any $y,z$
$$
\overline\Gamma^\epsilon_{zy} \big(f_x^{\epsilon;\tau}(y)\big) = f_x^{\epsilon;\tau}(z),
$$ 
this function belongs to $\mcD^\infty$ and its reconstruction satisfies 
$$
\overline\mcR^\epsilon \big(f_x^{\epsilon;\tau}\big)(y) = \overline\Pi^\epsilon_x(\tau)(y)
$$ 
regardless of the sign of $r(\tau)$. As an analogue, we use the following relations to define modelled distributions $\widetilde{f}_x^{\epsilon;\tau} : \bbR^d\to\widetilde{\bfT}$ representing $\nabla_h\overline{\Gamma}_{yx}^\epsilon(\tau)$, for any fixed $h\in H$. The dependence of $\widetilde{f}_x^{\epsilon;\tau}$ on $h$ is not emphasized in the notation for this function.
\begin{enumerate}
	\item We set
	\begin{align*}
	\widetilde{f}_x^{\epsilon;\Xi}(y) &\defeq \dot{\Xi}+h^\epsilon(y){\bf1},  \qquad  \widetilde{f}_x^{\epsilon;X^k} \defeq 0.
	\end{align*}
	In the first identity, the additional term $h^\epsilon(y){\bf1}$ is necessary to ensure that $\widetilde\mcR^{\epsilon,h}(\widetilde{f}_x^{\epsilon;\Xi}) = h^\epsilon$. (Recall that $\widetilde{\Pi}_x^{\epsilon,h}(\dot{\Xi}) = h^\epsilon-h^\epsilon(x)$.) The second identity reflects the fact that $\Gamma_{yx}(X^k)$ is deterministic.

	\item We require that some form of Leibniz rule holds true
	$$
	\widetilde{f}_x^{\epsilon;\tau\sigma} \defeq \widetilde{f}_x^{\epsilon;\tau} f_x^{\epsilon;\sigma} + f_x^{\epsilon;\tau}\widetilde{f}_x^{\epsilon;\sigma}.
	$$
	
	\item We use the notations
	$$
	\overline\mcJ^\epsilon(x)\tau \defeq \sum_{|k|<r(\mcI\tau)} \big(\partial^kK*\overline{\Pi}_x^\epsilon\tau\big)(x) \, \frac{X^k}{k!}
	$$
	and $Q_{<\alpha}:\widetilde{\bfT}\to\widetilde{\bfT}_{<\alpha}$ for the canonical projection, and set
	$$
	\widetilde{f}_x^{\epsilon;\mcI\tau}(y) \defeq \big(\mcI+ \overline\mcJ^\epsilon(y)\big)\widetilde{f}_x^{\epsilon;\tau}(y) - \overline{\Gamma}_{yx}^\epsilon \circ Q_{<r(\mcI\tau)} \big( \overline\mcJ^\epsilon(x)\widetilde{f}_x^{\epsilon;\tau}(x) \big).
	$$
	The truncation $Q_{<r(\mcI\tau)}$ is necessary because $\dot{\Xi}$ has a larger degree than $\Xi$.
	
	\item Last, we ask that
	$$
	\widetilde{f}_x^{\epsilon;\mcI_k\tau} \defeq \partial^k\widetilde{f}_x^{\epsilon;\mcI\tau},
	$$
	where the action of the operator $\partial^k\defeq\prod_{i=1}^d\partial_i^{k_i}$ on $X^\ell$ and planted trees $\mcI_\ell\tau$ is defined by
	$$
	\partial_iX^\ell \defeq {\bf1}_{\ell_i\ge1} \, \ell_i \, X^{\ell-e_i},  \qquad  \partial_i(\mcI_\ell\tau) = \mcI_{\ell+e_i}\tau.
	$$
\end{enumerate}	
Propositions 4.12 and 4.18 of \cite{HS} assert that one has the relation \eqref{HS:eq:realandformalderivative}.

\medskip

%%-----------------------%%
\noindent -- \textbf{Analytic step$\boldmath{.}$}
%%-----------------------%%
We cannot expect an $\epsilon$-uniform bound for $\widetilde{f}_x^{\epsilon;\tau}$ because the regularity of $h$ is overestimated. Indeed the degree $r(\dot{\Xi})=\theta>0$ is higher than the regularity exponent `$0$' of the space $L^2=B_{2,2}^0$. Instead we consider a truncation of $\widetilde{f}_x^{\epsilon;\tau}$ by removing some high order terms. For that purpose, for any $\tau\in\mcB$ denote by $\alpha_\tau$ the smallest degree of all non-polynomial symbols appearing in the smallest sector containing $\tau$, and set 
$$
\gamma_\tau\defeq\alpha_\tau + \frac{d}{2}.
$$

\begin{lem} \cite[Proposition 4.7]{HS} \label{LemmaProp4.7}
The function 
$$
\widetilde{F}^{\epsilon;\tau}_x \defeq Q_{<\gamma_\tau} \circ \widetilde{f}_x^{\epsilon;\tau}
$$ 
is a pointed modelled distribution in the class $\mcD_2^{\gamma_\tau,r(\tau)+d/2;x}$. As $\widetilde{f}_x^{\epsilon;\tau}$ itself, $\widetilde{F}^{\epsilon;\tau}_x$ depends implicitly on $h$. For each $i$ there exists a polynomial $P$ and a compact set $C\subset\bbR^d$ containing $x$ such that one has almost surely
$$
\sup_{\|h\|_H = 1}\tri \widetilde{F}^{\epsilon;\tau}_x \tri_{2,\gamma_\tau,r(\tau)+d/2;x}
\le
P\big(\tri\overline{\sf M}^\epsilon\tri_{\mcB_{i-1};C}\big).
$$
\end{lem}

\ssk

The function $\widetilde{F}^{\epsilon;\tau}_x$ is denoted $H^{x,h}_{\tau;n}$ in \cite{HS}. Alternatively, we can inductively define $\widetilde{F}^{\epsilon;\tau}_x$ from initial functions $\widetilde{F}^{\epsilon;\Xi}_x$ and $\widetilde{F}^{\epsilon;X^k}_x$ via two operations on pointed modelled distributions: tree products and integrations -- see Theorems 3.11 and 3.21 of \cite{HS}. The above lemma is proved via those inductive formulas. Lemma \ref{HS:LemmaProp4.6} below is partially involved in the inductive step when applying the integration map to $\widetilde{F}^{\epsilon;\tau}_x$ with $\gamma_\tau\le0$.

We return to the $h$-uniform estimate of $\nabla_h\overline{\Pi}_x^\epsilon(\tau_i)$. If we are in a case where $\gamma_{\tau_i}>0$ the $\widetilde\mcR^{\epsilon,h}$-reconstruction of $\widetilde{F}^{\epsilon;\tau}_x\in \mcD_2^{\gamma_\tau,r(\tau)+d/2;x}$ is well-defined and satisfies
$$
\nabla_h\overline{\Pi}_x^\epsilon(\tau_i) = \widetilde\mcR^{\epsilon,h}\big(\widetilde{f}_x^{\epsilon;\tau_i}\big) = \widetilde\mcR^{\epsilon,h} \big(\widetilde{F}^{\epsilon;\tau_i}_x \big),
$$ 
so an $(h,\epsilon)$-uniform estimate on a moment of $\big(\nabla_h\overline{\Pi}_x^\epsilon(\tau)\big)(\varphi_x^\lambda)$ follows by induction from Lemma \ref{LemmaProp4.7} and the refined reconstruction theorem given in Theorem \ref{HS:thm:pointreconst}. We cannot use Theorem \ref{HS:thm:pointreconst} if $\gamma_{\tau_i}\le0$ since there is not a unique reconstruction in that case. However the following lemma shows that $ \widetilde\mcR^{\epsilon,h}\big(\widetilde{f}_x^{\epsilon;\tau_i}\big)$ satisfies a good estimate as a candidate for a reconstruction of $\widetilde{F}_x^{\epsilon;\tau_i}$.

\medskip

\begin{lem}\label{HS:LemmaProp4.6}
\cite[Lemma 4.6]{HS}
Even if $\gamma_{\tau_i}\le0$ there exists a polynomial $P$ and a compact set $C\subset\bbR^d$ containing $x$ such that one has almost surely
$$
\sup_{\|h\|_H=1} m\Big( \widetilde{F}^{\epsilon;\tau_i}_x , \widetilde\mcR^{\epsilon,h}(\widetilde{f}_x^{\epsilon;\tau_i}) ; x \Big) \lesssim P\big(\tri\overline{\sf M}^\epsilon\tri_{\mcB_{i-1};C}\big).
$$
\end{lem}

\ssk

When $\tau_i=\Xi$ the above lemma can be checked directly from some elementary deterministic estimates on $h^\varepsilon$. Otherwise $\tau$ should be of the form $\tau=\Xi\sigma$ with some $\sigma$ such that $\alpha_\sigma+r(\Xi)+d/2>0$. This is the content of Lemma 4.4 of \cite{HS}. This condition makes it possible to prove the desired estimate as a result of some kind of `Young multiplication' between $B_{2,\infty}^{r(\Xi)+d/2}$ and $B_{\infty,\infty}^{\alpha_\sigma}$, as in Theorem 3.12 of \cite{BL22}. The bound on $\widetilde{F}^{\epsilon;\tau}_x$ from Lemma \ref{LemmaProp4.7} can then be transferred to $\big(\widetilde\mcR^{\epsilon,h}(\widetilde{f}_x^{\epsilon;\tau_i})\big)(\varphi_x^\lambda) = (\nabla_h\overline{\Pi}_x^\epsilon\tau)(\varphi_x^\lambda)$, so the result follows by induction.

\medskip

%----------------------------------------------------------------------------------%
\subsection{Bailleul \& Hoshino \cite{RandomModel} and regularity-integrability structures$\boldmath{.}$ \hspace{0.1cm}}
\label{SectionBH}
%----------------------------------------------------------------------------------%

We describe in this section the approach of \cite{RandomModel} for proving (an improved version of) Theorem \ref{HS:thm:main}.

%%-----------------------------------------------------------------------%%
\subsubsection{The abstract Malliavin derivative and an example$\boldmath{.}$ \hspace{0.1cm}}
\label{SubsectionBH1}
%%-----------------------------------------------------------------------%%

We use the extended model space $\widetilde{\bfT}$ introduced in Section \ref{SubsectionHS1}, including the symbol $\dot{\Xi}$. Unlike \cite{HS} we do not assign here a positive degree to $\dot \Xi$. Rather, noting that $H=L^2(\bbR^d)$ embeds into $B_{\infty,\infty}^{-d/2}$ we assign $\dot{\Xi}$ the same degree $\alpha_0<-d/2$ as $\Xi$. With that choice, it turns out that one can build a setting where we have the intertwining relation
\begin{align}\label{BH:eq:realandformalderivative}
\nabla_h\overline{\Pi}_x^\epsilon(\tau) = \widehat{\Pi}_x^{\epsilon,h}(D\tau)
\end{align}
for all $\tau\in\mcB$, for some model $\widehat{\sf M}^{\epsilon,h}$ introduced below in Section \ref{SubsectionBH3} as $\widehat{\sf M}^{\epsilon,h;\infty}$. This identity plays in \cite{RandomModel} a key role in proving some estimates on $\nabla_h(\overline\Pi_x^\epsilon\tau)(\varphi_x^\lambda)$ which, together with the spectral gap inequality \eqref{HS:useofpoincare}, allow for an inductive proof of Theorem \ref{HS:thm:main}. We explain in this paragraph the mechanics on the example of the symbol $\begin{tikzpicture}[scale=0.3,baseline=0.05cm]
	\node at (0,0)  [noise] (1) {};
	\node at (0,1)  [noise] (2) {};
	\draw[K] (1) to (2);
\end{tikzpicture}=\Xi\mcI(\Xi)$
in dimension $d=3$. 
(The following estimates are locally uniform over $x$, but for simplicity, we describe them as if they were global over $x$.) Since $r(\begin{tikzpicture}[scale=0.3,baseline=0.05cm]
	\node at (0,0)  [noise] (1) {};
	\node at (0,1)  [noise] (2) {};
	\draw[K] (1) to (2);
\end{tikzpicture}) =-1-$, the desired estimate is the $\{\Vert h\Vert_H=1\}$-uniform bound
\begin{align}\label{BH:ex:goal}
\sup_x \big| \big(\nabla_h\overline{\Pi}_x^\epsilon( \begin{tikzpicture}[scale=0.3,baseline=0.05cm]
	\node at (0,0)  [noise] (1) {};
	\node at (0,1)  [noise] (2) {};
	\draw[K] (1) to (2);
\end{tikzpicture})\big)(\varphi_x^\lambda) \big| \lesssim \lambda^{-1-}.
\end{align}
We start from the identity
\begin{align*}
\nabla_h\overline{\Pi}_x^\epsilon\tau = \widehat{\Pi}_x^{\epsilon,h}(D\tau)
&= \big\{K*\xi^\epsilon - K*\xi^\epsilon(x)\big\} h^\epsilon + \big\{K*h^\epsilon - K*h^\epsilon(x)\big\}\xi^\epsilon   \\
&\eqdef A_x^\epsilon + B_x^\epsilon.
\end{align*}
One can use a Besov-type reconstruction argument to bound $A_x^\epsilon$ -- see Theorem 3.2 of \cite{BL22} for that type of reconstruction theorem and Theorem 4.4 of \cite{CZ} for the initial version of the reconstruction theorem for coherent germs. Indeed the germ $f^\epsilon_x\defeq(K*\xi^\epsilon)(x)h^\epsilon$ satisfies a $B_{2,\infty}^{1/2-}$-type coherence property
\begin{align*}
\big\| (f^\epsilon_{x+y} - f^\epsilon_x)(\varphi_x^\lambda) \big\|_{L_x^2} = \big\|\big(K*\xi^\epsilon(x+y) - K*\xi^\epsilon(x)\big)h^\epsilon(\varphi_x^\lambda)\big\|_{L_x^2}
\lesssim|y|^{1/2-},
\end{align*}
since the family $\{K*\xi^\epsilon\}_{0<\epsilon\leq 1}$ is bounded in $B_{\infty,\infty}^{1/2-}$ and the family $\{h^\epsilon\}_{0<\epsilon\leq 1}$ is bounded in $L^2$. Since the reconstruction of the germ $f^\epsilon=\{f^\epsilon_x\}_x$ is given by $\mcR f^\epsilon = (K*\xi^\epsilon)h^\epsilon$, we have the estimate
$$
\big\| A_x^\epsilon(\varphi_x^\lambda) \big\|_{L_x^2} = \big\| \big(\mcR f^\epsilon - f^\epsilon_x\big)(\varphi_x^\lambda) \big\|_{L_x^2} \lesssim \lambda^{1/2-}
$$
from the characterization of the reconstruction of a germ. Although this `$B_{2,\infty}^{1/2-}$-type' estimate is different from \eqref{BH:ex:goal}, by noting that $B_{2,\infty}^{1/2-}$ is embedded into $B_{\infty,\infty}^{-1-}$, we expect that a `$B_{\infty,\infty}^{-1-}$-type' estimate
\begin{align}\label{BH:ex:goal1}
\sup_x| A_x^\epsilon(\varphi_x^\lambda) | \lesssim \lambda^{-1-}
\end{align}
also holds as a consequence of a similar embedding argument. This heuristic argument is justified by Lemma 12 of \cite{RandomModel}. A similar reasoning does not work to obtain the desired estimate for $B_x^\epsilon$. Indeed the germ $g^\epsilon_x\defeq (K*h^\epsilon)(x)\xi^\epsilon$ has only a weaker coherence property
\begin{align*}
\big\| (g^\epsilon_{x+y} - g^\epsilon_x)(\varphi_x^\lambda) \big\|_{L_x^2} = \big\|\big(K*h^\epsilon(x+y)-K*h^\epsilon(x)\big)\xi^\epsilon(\varphi_x^\lambda)\big\|_{L_x^2}
\lesssim|y|\lambda^{-3/2-},
\end{align*}
because the difference $K*h^\epsilon(x+y)-K*h^\epsilon(x)$ cannot produce a superlinear estimate of $y$ even if $\{K*h^\epsilon\}_{0<\epsilon\leq 1}$ is bounded in $B_{2,\infty}^2$.
Instead, the germ $g^{\epsilon+}$ defined by
$$
g^{\epsilon+}_x\defeq \bigg\{(K*h^\epsilon)(x)+\sum_{i=1}^d(\cdot-x)_i \, \partial_iK*h^\epsilon(x)\bigg\}\xi^\epsilon,
$$
has the desired coherence property
$$
\big\| \big(g^{\epsilon+}_{x+y} - g^{\epsilon+}_x\big)(\varphi_x^\lambda) \big\|_{L_x^2} \lesssim |y|^2\lambda^{-3/2-}.
$$
Since the reconstruction of $g^{\epsilon+}$ is given by $\mcR g^{\epsilon+} = (K*h^\epsilon)\xi^\epsilon$, we obtain a $B_{2,\infty}^{1/2-}$-type estimate
$$
\| B^{\epsilon+}_x(\varphi_x^\lambda) \|_{L_x^2} \lesssim \lambda^{1/2-}
$$
for $B^{\epsilon+}_x \defeq \big\{ K*h^\epsilon(\cdot)-K*h^\epsilon(x)-\sum_i(\cdot-x)_i \, \partial_iK*h^\epsilon(x) \big\} \xi^\epsilon(\cdot)$. To eventually get an estimate on $B_x^\epsilon(\varphi_x^\lambda)$ we need some control on the terms
$$
C_{x,i}^{\epsilon}\defeq(\cdot-x)_i \, \partial_iK*h^\epsilon(x)\xi^\epsilon(\cdot).
$$
Note that $(\cdot-x)_i\xi^\epsilon$ is equal to the application of $\overline{\Pi}_x^\epsilon$ to the symbol $X_i\Xi$ with degree $-1/2-$, and the family $\{\partial_iK*h^\epsilon\}_{0<\epsilon\leq 1}$ is bounded in $B_{2,\infty}^1\subset L^6$. Thus we have a `$B_{6,\infty}^{-1/2-}$-type' estimate
$$
\|C_{x,i}^\epsilon(\varphi_x^\lambda)\|_{L_x^6} = \big\| (\partial_iK*h^\epsilon)(x)\overline{\Pi}_x^\epsilon(X_i\Xi)(\varphi_x^\lambda) \big\|_{L_x^6}
\lesssim \lambda^{-1/2-}.
$$
Since $B_{6,\infty}^{-1/2-}$ is embedded into $B_{\infty,\infty}^{-1-}$, we expect that a $B_{\infty,\infty}^{-1-}$-type estimate similar to \eqref{BH:ex:goal1} also holds for $B_x^\varepsilon=B_x^{\varepsilon+}+\sum_iC_{x,i}^\varepsilon$ by an embedding argument similar to the argument used for the $A_x^\epsilon$ term.

\ssk

This example illustrate an important insight. It is useful to adjust the regularity and integrability exponents of $h$ in different computations: we considered $h\in B_{\infty,\infty}^{-3/2}$ to show the identity \eqref{BH:eq:realandformalderivative}, we used that $h\in L^2$ in the estimates on $A_x^\epsilon$ and $B_x^{\epsilon+}$, and we used that $h\in B_{6,\infty}^{-1}$ in the estimate on $C_{x,i}^\epsilon$. To place the above argument in a more general context we consider an extension of the notion of regularity structure including some integrability exponents and a hierarchy based on the numerology of the classical Besov embeddings.

\medskip

%%-----------------------%%
\subsubsection{Regularity-integrability structures$\boldmath{.}$ \hspace{0.1cm}}
\label{SubsectionBH2}
%%-----------------------%%

We use the bold letters $\bfa,\bfb,\bfc$ to denote some generic elements of $\bbR\times[1,\infty]$. We represent each component of $\bfa\in\bbR\times[1,\infty]$ by 
$$
\bfa = \big(r(\bfa),i(\bfa)\big),
$$
where the letters ``$r$" and ``$i$" stand for ``regularity" and ``integrability", respectively. We define a partial order $\preceq$ and a strict partial order $\prec$ on the set $\bbR\times[1,\infty]$ by setting
\begin{align*}
\bfb\preceq\bfa&\quad\overset{\text{def}}{\Leftrightarrow}\quad r(\bfb)\le r(\bfa),\ i(\bfb)\ge i(\bfa),\\
\bfb\prec\bfa&\quad\overset{\text{def}}{\Leftrightarrow}\quad r(\bfb)<r(\bfa),\ i(\bfb)\ge i(\bfa).
\end{align*}
Note that $i(\bfb)$ may be equal to $i(\bfa)$ even for the latter case. For any $\bfb\preceq\bfa$, we define the element $\bfa\ominus\bfb\in\bbR\times[1,\infty]$ by
$$
\bfa\ominus\bfb\defeq\Bigg(r(\bfa)-r(\bfb),\frac1{\frac1{i(\bfa)}-\frac1{i(\bfb)}}\Bigg).
$$
We have in particular the relations $r(\bfa)=r(\bfa\ominus\bfb)+r(\bfb)$ and $1/{i(\bfa)}=1/{i(\bfa\ominus\bfb)}+1/{i(\bfb)}$.

\medskip

\begin{definition*}
A {\bf regularity-integrability structure} $\scT=(\bfA,\bfT,\bfG)$ consists of 
\begin{itemize} \setlength{\itemsep}{0.1cm}
\item[(1)]
$\bfA$ a subset of $\bbR\times[1,\infty]$ such that, for every $\bfc\in\bbR\times[1,\infty]$, the subset $\{\bfa\in \bfA\, ;\, \bfa\prec\bfc\}$ is finite.
\item[(2)]
$\bfT=\bigoplus_{\bfa\in \bfA}\bfT_\bfa$ an algebraic sum of Banach spaces $(\bfT_\bfa,\|\cdot\|_\bfa)$ indexed by $\bf A$.
\item[(3)]
$\bfG$ a group of continuous linear operators on $\bfT$ such that, for any $\Gamma\in \bfG$ and $\bfa\in\bfA$,
$$
(\Gamma-\id)\bfT_\bfa\subset \bfT_{\prec\bfa}\defeq\bigoplus_{\bfb\in \bfA,\, \bfb\prec\bfa}\bfT_\bfa.
$$
\end{itemize}
The biggest number $\alpha_0\in\bbR$ satisfying $(\alpha_0,\infty)\preceq\bfa$ for any $\bfa\in\bfA$ is called the \emph{regularity of} $\scT$. For any $\bfa\in\bfA$ we denote by $P_\bfa:\bfT\to\bfT_\bfa$ the canonical projection and write
$$
\|\tau\|_\bfa \defeq \|P_\bfa\tau\|_\bfa  \qquad  (\tau\in \bfT).
$$
\end{definition*}

\medskip

The notions of models and modelled distributions are defined accordingly.

\medskip

\begin{definition*} 
Let $\scT=(\bfA,\bfT,\bfG)$ be a regularity-integrability structure of regularity $\alpha_0$. Fix a positive integer $r>|\alpha_0|$. A {\bf model} ${\sf M}=(\Pi,\Gamma)$ for $\scT$ on $\bbR^d$ consists of two families of continuous linear operators
$$
\Pi = \big\{ \Pi_x:\bfT\to \mcD'(\bbR^d) \big\}_{x\in\bbR^d},  \qquad   \Gamma = \{\Gamma_{yx}\}_{x,y\in\bbR^d}\subset \bfG
$$
satisfying the following properties.
\begin{itemize} \setlength{\itemsep}{0.1cm}
\item[(1)]
We have $\Pi_x\Gamma_{xy}=\Pi_y$, $\Gamma_{xx}=\id$, $\Gamma_{xy}\Gamma_{yz}=\Gamma_{xz}$ for any $x,y,z\in\bbR^d$.
\item[(2)]
For any $\bfc\in\bbR\times[1,\infty]$ and any compact subset $C$ of $\bbR^d$
\begin{align*}
\|\Pi\|_{\bfc;C}&\defeq
\max_{\bfa\in\bfA,\,\bfa\prec\bfc}\sup_{\varphi\in\mcB_r}\sup_{\lambda\in(0,1]}
\lambda^{-r(\bfa)}
\bigg\|\sup_{\tau\in\bfT_\bfa\setminus\{0\}}
\frac{|(\Pi_x\tau)(\varphi_x^\lambda)|}{\|\tau\|_\bfa}\bigg\|_{L_x^{i(\bfa)}(C)}<\infty,
\\
\|\Gamma\|_{\bfc;C}&\defeq
\max_{\bfa,\bfb\in\bfA,\,\bfb\prec\bfa\prec\bfc}
\sup_{0<\|y\|\le1}\frac1{\|y\|^{r(\bfa\ominus\bfb)}}
\bigg\|
\sup_{\tau\in\bfT_\bfa\setminus\{0\}}
\frac{\|\Gamma_{(x+y)x}(\tau)\|_\bfb}
{\|\tau\|_\bfa}
\bigg\|_{L_x^{i(\bfa\ominus\bfb)}(C)}
<\infty.
\end{align*}
\end{itemize}
We define $\tri {\sf M}\tri_{\bfc;C}$ and $\tri {\sf M}^{(1)} \,,\, {\sf M}^{(2)}\tri_{\bfc;C}$ as in Definition \ref{HS:def:model}.
\end{definition*}

\medskip

Although we employ here some local estimates, in line with the previous sections, in the papers \cite{Hos23} and \cite{RandomModel} where the regularity-integrability structures are introduced, we use some global estimates involving some heat kernels and weight functions. This difference is merely technical and not a serious matter. The same remark applies to the following definition of a modelled distribution in a regularity-integrability setting.

\medskip

\begin{definition*}
Let $\scT=(\bfA,\bfT,\bfG)$ be a regularity-integrability structure and let ${\sf M} = (\Pi,\Gamma)$ be a model for $\scT$ on $\bbR^d$. For any $\bfc\in\bbR\times[1,\infty]$ we define $\mcD^\bfc=\mcD^\bfc(\Gamma)$ as the space of all functions $f:\bbR^d\to\bfT_{\prec\bfc}$ for which
\begin{align*}
\rest{f}_{\bfc;C}&:=
\max_{\bfa\in\bfA,\,\bfa\prec\bfc}\big\|\|f(x)\|_\bfa\big\|_{L_x^{i(\bfc\ominus\bfa)}(C)}<\infty,
\\
\|f\|_{\bfc;C}&:=
\max_{\bfa\in\bfA,\,\bfa\prec\bfc}
\sup_{0<\|y\|\le1}\frac{\big\| \big\| f(x+y)-\Gamma_{(x+y)x}(f(x)) \big\|_\bfa \big\|_{L_x^{i(\bfc\ominus\bfa)}(C)}}{\|y\|^{r(\bfc\ominus\bfa)}}<\infty
\end{align*}
for any compact subset $C\subset\bbR^d$. We define $\tri {\sf M}\tri_{\bfc;C}$ and $\tri {\sf M}^{(1)} \,,\, {\sf M}^{(2)}\tri_{\bfc;C}$ as in Definition \ref{HS:def:MD}.
\end{definition*}

\medskip

For any $f\in\mcD^\bfc$ the germ $\big\{f_x\defeq\Pi_x(f(x))\big\}_{x\in\bbR^d}$ satisfies the `coherence property' in the sense of \cite{BL22}. The following reconstruction theorem is thus a consequence of Theorem 3.2 of \cite{BL22}. (See Theorem 4.1 of \cite{Hos23} for the proof in a slightly different setting.)

\medskip

\begin{thm} \label{BH:thm:reconsRIS}
Let $\scT=(\bfA,\bfT,\bfG)$ be a regularity-integrability structure and let $M=(\Pi,\Gamma)$ be a model for $\scT$ on $\bbR^d$. For any $\bfc\in(0,\infty)\times[1,\infty]$ there exists a unique continuous linear map $\mcR^{\sf M}:\mcD^\bfc(\Gamma)\to\mcD'(\bbR^d)$ satisfying
$$
\sup_{\varphi\in\mcB_r}\sup_{\lambda\in(0,1]}\lambda^{-r(\bfc)} \big\|\big(\mcR^{\sf M} f-\Pi_xf(x)\big)(\varphi_x^\lambda)\big\|_{L_x^{i(\bfc)}(C)}\lesssim\|\Pi\|_{\bfc;\overline{C}}\|f\|_{\bfc;\overline{C}}.
$$
Moreover the mapping $({\sf M},f)\mapsto\mcR^{\sf M} f$ is locally Lipschitz continuous with respect to the quasi-metrics $\tri {\sf M}^{(1)} \,,\, {\sf M}^{(2)}\tri_{\bfc;\overline{C}}$ and $\tri f^{(1)},f^{(2)}\tri_{\bfc;\overline{C}}$.
\end{thm}

\medskip

%%----------------------------------------------------------------------------------%%
\subsubsection{The convergence result in a regularity-integrability setting$\boldmath{.}$ \hspace{0.1cm}}
\label{SubsectionBH3}
%%----------------------------------------------------------------------------------%%

Recall $\dot{\mcB}$ consists of the set of symbols obtained by replacing one argument $\Xi$ in $\tau\in\mcB$ with a $\dot{\Xi}$. For any $p\in[2,\infty]$ we define the extended degree map $r_p:\mcB\cup\dot{\mcB}\to\bbR$ by the same rules as \eqref{Basics:defofr} with the additional one
$$
r_p(\dot{\Xi})=\alpha_0+\frac{d}p.
$$
The value of $r_p(\dot{\Xi})$ comes from the numerology in the Besov embedding $L^2\subset B_{2,\infty}^{\alpha_0+d/2}\subset B_{p,\infty}^{\alpha_0+d/p}\subset B_{\infty,\infty}^{\alpha_0}$.
 The map $i_p:\mcB\cup\dot{\mcB}\to[2,\infty]$ is defined by
$$
i_p(\tau)=
\begin{cases}
\infty&(\tau\in\mcB),\\
p&(\tau\in\dot{\mcB}).
\end{cases}
$$
%Roughly speaking, our aim is to show the $B_{i_p(\tau),\infty}^{r_p(\tau)}$-type estimate for the application of the model to each $\tau\in\mcB\cup\dot{\mcB}$. 
The set $\big\{\tau\in\mcB\, ;\, (r_p(\tau),i_p(\tau))\prec\bfc\big\}$ is finite for any $\bfc\in\bbR\times[1,\infty]$ and we can construct a regularity-integrability structure 
$$
\widehat\scT_p = \big( \widehat\bfA \,,\, \widehat\bfT \,,\, \widehat\bfG \big)
$$ 
by setting
\begin{itemize} \setlength{\itemsep}{0.1cm}
\item[(1)]
$\widehat\bfA = \big\{ (r_p(\tau),i_p(\tau))\,;\,\tau\in\mcB\cup\dot{\mcB} \big\}$,
\item[(2)]
$\widehat\bfT = \bigoplus_{\bfa\in \bfA}\bfT_\bfa$, where $\widehat\bfT_\bfa = \spa\big\{ \tau\in\mcB\cup\dot\mcB\,;\, (r_p(\tau),i_p(\tau))=\bfa \big\}$,
\item[(3)]
$\widehat\bfG$ is the set of linear maps $\Gamma : \widehat\bfT\to\widehat\bfT$ satisfying for each $\tau\in\mcB\cup\dot{\mcB}$
$$
(\Gamma-\id)\tau \in \spa\Big\{ \sigma\in\mcB\cup\dot{\mcB} \,;\, (r_p(\sigma),i_p(\sigma))\prec(r_p(\tau),i_p(\tau)) \Big\}.
$$
\end{itemize}
(The space $\widehat\scT_p$ is denoted by $\scW_p$ in \cite{RandomModel}.) There is a unique multiplicative map ${\bf \Pi}^{\epsilon,h}$ on $\widehat\bfT = \spa(\mcB\cup\dot\mcB)$ such that $({\bf\Pi}^{\epsilon,h}X^k)(x)=x^k, \; {\bf\Pi}^{\epsilon,h}(\Xi)=\xi^\epsilon, \; {\bf\Pi}^{\epsilon,h}(\dot{\Xi})=h^\epsilon$ and 
$$
{\bf\Pi}^{\epsilon,h}(\mcI_k\tau) = \partial^kK*{\bf\Pi}^{\epsilon,h}\tau \qquad (\tau\in \bfT).
$$
This map induces a unique admissible model ${\sf M}^{\epsilon,h;p}$ on $\widehat\scT_p$ which depends on $p$ via the degree map $r_p$. This is for instance visible in the identity
$$
\Pi_x^{\epsilon,h;p}(\mcI\tau) = K*(\Pi_x^{\epsilon,h;p}\tau) - \sum_{|\ell|<r_p(\mcI_k(\tau))}\frac{(\cdot-x)^\ell}{\ell!} \, \big(\partial^\ell K*\Pi_x^{\epsilon,h;p}\tau\big)(x).
$$
We keep from that expression the fact that the order of the Taylor expansion varies with $p$.

As in Section \ref{SubsectionBasicsBPHZ}, the class of preparation maps provides a systematic way of constructing some admissible models. One can in particular use the same type of map
\begin{equation} \label{EqTildeR}
\widehat{R}_\ell = (\ell\otimes\id)\widehat\Delta_r^-
\end{equation}
as in Section \ref{SubsectionBasicsBPHZ}, with $\widehat\Delta_r^-$ in place of $\Delta_r^-$, where $\widehat\Delta_r^-$ is an extension of $\Delta_r^-$ which cuts any given tree into a subtree with negative degree which contains the root of $\tau$ and does not contain $\dot{\Xi}$, and the remaining graph. The model $\widehat{\sf M}^{\epsilon,h;p}$ on $\widehat{\scT}_p$ associated with the preparation map $\widehat{R}_\ell$ is also called the BPHZ model. Its convergence to a limit model entails the convergence of the BPHZ model $\overline{\sf M}^\epsilon$ of Section \ref{SubsectionBasicsBPHZ} to a limit. Practically, the main result of \cite{RandomModel} implies that if $\Xi$ is the only symbol of $\mcB$ with degree less than or equal to $-d/2$ then for any $\bfc\in\bbR\times[1,\infty]$, $p\in[2,\infty]$, $q\ge1$ and any compact set $C$ of $\bbR^d$ we have
\begin{align} \label{thm:main1}
\lim_{\epsilon_1,\epsilon_2\to 0}
\bbE\bigg[\sup_{\|h\|_H = 1}
\tri\widehat{\sf M}^{\epsilon_1,h;p},\widehat{\sf M}^{\epsilon_2,h;p}\tri_{\bfc;C}^q
\bigg]=0.
\end{align}
Theorem 8 in \cite{RandomModel} actually provides the same result for a class of models associated with a family of preparation maps that contains the BPHZ preparation map \eqref{EqTildeR} as a particular case.

\ssk

%%--------------------------------------------------------------%%
\subsubsection{A sketch of the proof of \eqref{thm:main1}$\boldmath{.}$ \hspace{0.1cm}}
\label{SubsectionBH4}
%%--------------------------------------------------------------%%

We use the same ordering of $\mcB$ as in Section \ref{SubsectionHS3}
$$
\mcB\backslash\{X^k\}_{k\in\bbN^d} = \big\{\tau_1 \preceq \tau_2 \preceq \cdots\big\}.
$$ 
For each $i$ we define $\dot{\mcB}_i$ as the set of all symbols obtained by replacing one argument $\Xi$ in $\tau\in\mcB_i$ with $\dot{\Xi}$. Similarly to Section \ref{SubsectionHS3}, we only outline here the proof of the uniform bounds. A simple modification provides the proof of the convergence as all the identities below have some locally Lipschitz counterparts. This makes the proof of convergence easier than in \cite{HS} as we do not need to introduce any new notion of modelled distributions measured in some negative Sobolev norms.

\ssk

For any finite subset $\mcA\subset\mcB\cup\dot{\mcB}$ and any $p\in[2,\infty]$ we denote by $\text{\bf bd}(\mcA,p)$ the statement that for any $q\ge1$ and any compact set $C\subset\bbR^d$ we have
$$
\sup_{0<\epsilon\leq 1}
\bbE\bigg[
\sup_{\|h\|_{H}=1}
\|\widehat{\sf M}^{\epsilon,h;p}\|_{\mcA;C}^q
\bigg]<\infty.
$$
We also denote by $\text{\bf bd}(\mcA)$ the statement that $\text{\bf bd}(\mcA,p)$ holds for all $p\in[2,\infty]$. The flow of the proof is explained as follows.
\begin{itemize} \setlength{\itemsep}{0.1cm}
\item[(1)]
$\text{\bf bd}(\dot{\mcB}_1)$ follows from some elementary deterministic estimates on $h^\epsilon$.
\item[(2)]
We show $\text{\bf bd}(\mcB_i)$ and $\text{\bf bd}(\dot{\mcB}_i)$ for any $i$ in the following steps.
\begin{itemize}
\item[(a)]
{\bf Probabilistic step}:
$\text{\bf bd}(\mcB_{i-1}\cup\dot{\mcB}_i) \Longrightarrow \text{\bf bd}(\mcB_i)$.
\item[(b)]
{\bf Analytic step}:
$\text{\bf bd}(\mcB_i\cup\dot{\mcB}_i,2) \Longrightarrow \text{\bf bd}(\dot{\mcB}_{i+1},2)$.
\item[(c)]
{\bf Algebraic step}:
$\text{\bf bd}(\mcB_i\cup\dot{\mcB}_i)$ \& $\text{\bf bd}(\dot{\mcB}_{i+1},2) \Longrightarrow \text{\bf bd}(\dot{\mcB}_{i+1})$.
\end{itemize}
\end{itemize}
The following diagram illustrates the mechanics; the assertions are shown in the order indicated by the solid arrows.
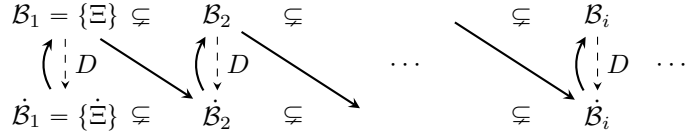
\begin{center}
\begin{tikzpicture}[auto]
\node (a1) at (0,0) {$\mcB_1=\{\Xi\}$}; 
\node (a2) at (2,0) {$\mcB_2$}; 
\node (a3) at (4,0) {}; 
\node (a99) at (5,0) {}; 
\node (a100) at (7,0) {$\mcB_i$}; 
\node (a101) at (8,0) {}; 
\node (b1) at (0,-1.3) {$\dot{\mcB}_1=\{\dot{\Xi}\}$}; 
\node (b2) at (2,-1.3) {$\dot{\mcB}_2$}; 
\node (b3) at (4,-1.3) {}; 
\node (b99) at (5,-1.3) {}; 
\node (b100) at (7,-1.3) {$\dot{\mcB}_i$}; 
\node (b101) at (8,-1.3) {}; 
\node (c1) at ($0.5*(a1)+0.5*(a2)$) {$\subsetneq$};
\node (d1) at ($0.5*(b1)+0.5*(b2)$) {$\subsetneq$};
\node (c2) at ($0.5*(a3)+0.5*(a2)$) {$\subsetneq$};
\node (d2) at ($0.5*(b3)+0.5*(b2)$) {$\subsetneq$};
\node (c99) at ($0.5*(a99)+0.5*(a100)$) {$\subsetneq$};
\node (d99) at ($0.5*(b99)+0.5*(b100)$) {$\subsetneq$};
\node (center) at ($0.25*(a3)+0.25*(a99)+0.25*(b3)+0.25*(b99)$) {$\cdots$};
\node (center2) at ($0.5*(a101)+0.5*(b101)$) {$\cdots$};
\draw[->, dashed] (a1) to node {$D$} (b1);
\draw[->, out=120, in=240, thick] (b1) to (a1);
\draw[->, thick] (a1) to (b2);
\draw[->, dashed] (a2) to node {$D$} (b2);
\draw[->, out=120, in=240, thick] (b2) to (a2);
\draw[->, thick] (a2) to (b3);
\draw[->, thick] (a99) to (b100);
\draw[->, dashed] (a100) to node {$D$} (b100);
\draw[->, out=120, in=240, thick] (b100) to (a100);
\end{tikzpicture}
\end{center}

\ssk

For readers familiar with the work \cite{LOTT} of Linares, Otto, Tempelmayr \& Tsatsoulis one can make the following parallel between the above reasoning and the mechanics of \cite{LOTT}. The step 2(a) is common to both approaches: this is where we use crucially the spectral gap ingredient. The Step 2(b) is somewhat the equivalent of their {\it algebraic} \& {\it three point} arguments. The Step 2(c) is the equivalent of their {\it Reconstruction III} step together with their {\it Averaging} step.

\medskip

%%-----------------------%%
\subsubsection{Probabilistic step: $\text{\bf bd}(\mcB_{i-1}\cup\dot{\mcB}_i) \Longrightarrow \text{\bf bd}(\mcB_i)$$\boldmath{.}$ \hspace{0.1cm}}
\label{SubsectionBH4.1}
%%-----------------------%%

Recall $\mcB_i=\mcB_{i-1}\cup\{\tau_i\}$. We show the estimate on $\widehat{\Pi}_x^{n,h;p}(\tau_i)$ assuming $\text{\bf bd}(\mcB_{i-1}\cup\dot{\mcB}_i)$. If $r_\infty(\tau_i)>0$ the estimate follows from the reconstruction theorem, in the same way as in the proof of Lemma \ref{HS:lem:reconst}. In that case the reasoning is purely deterministic. If $r_\infty(\tau_i)\le0$ we use the inequality \eqref{HS:useofpoincare} derived from Poincar\'e inequality, so the reasoning is probabilistic. The estimate of the expectation term follows from the property of BPHZ models. It is precisely at that point that the translation invariance of both the kernel and the law of the noise play a crucial role. For the derivative term we use the identity $\nabla_h\overline{\Pi}_x^\epsilon(\tau_i) = \widehat{\Pi}_x^{\epsilon,h;\infty}(D\tau_i)$ to deduce the estimate from the assumption $\text{\bf bd}(\dot{\mcB}_i)$ and the induction assumption. This is the content of the following statement.

\medskip

\begin{lem} \label{lem:BH4.1.2}
\cite[Lemma 10]{RandomModel}
If $r_\infty(\tau_i)\leq 0$ there exists a polynomial $P$ such that for any $q\ge1$ and any compact set $C\subset\bbR^d$ we have
$$
\bbE\big[\|\overline{\Pi}^\epsilon\|_{\tau_i;C}^q\big]^{\frac{1}{q}}
\le \bbE\bigg[P\bigg(\sup_{\|h\|_H = 1}\|\widehat{\Pi}^{\epsilon,h;\infty}\|_{\mcB_{i-1}\cup\dot{\mcB}_i;\overline{C}}\bigg)\bigg].
$$
\end{lem}

%\ssk

%%-----------------------%%
\subsubsection{Analytic step: $\text{\bf bd}(\mcB_i\cup\dot{\mcB}_i,2) \Longrightarrow \text{\bf bd}(\dot{\mcB}_{i+1},2)$$\boldmath{.}$ \hspace{0.1cm}}
\label{SubsectionBH4.2}
%%-----------------------%%

We show the estimate on $\widehat{\Pi}_x^{\epsilon,h;2}(\tau)$ for any $\tau\in\dot{\mcB}_{i+1}\setminus\dot{\mcB}_i$. Since $\tau$ is not equal to $\dot{\Xi}$ we have
$$
r_2(\tau)=r_\infty(\tau)+\frac{d}2>0
$$
from the assumption \eqref{HS:thm:main:asmp}. We can obtain the estimate on $\widetilde{\Pi}_x^{\epsilon,h;2}(\tau)$ by an argument similar to the argument used in the proof of Lemma \ref{HS:lem:reconst}, applying the regularity-integrability version of the reconstruction theorem, Theorem \ref{BH:thm:reconsRIS}, to a well-chosen modelled distribution.

\medskip

\begin{lemma*} \label{lem:BH4.2}
\cite[Lemma 15]{RandomModel}
There exists a polynomial $P$ such that for any $h\in H$ and any compact set $C\subset\bbR^d$ we have
$$
\|\widehat{\Pi}^{\epsilon,h;2}\|_{\dot{\mcB}_{i+1};C} \le P\big(\|\widehat{\Pi}^{\epsilon,h;2}\|_{\mcB_i\cup\dot{\mcB}_i;\overline{C}}\big).
$$
\end{lemma*}

\ssk

%%-----------------------%%
\subsubsection{Algebraic step: $\text{\bf bd}(\mcB_i\cup\dot{\mcB}_i)$ \& $\text{\bf bd}(\dot{\mcB}_{i+1},2) \Longrightarrow \text{\bf bd}(\dot{\mcB}_{i+1})$$\boldmath{.}$ \hspace{0.1cm}}
\label{SubsectionBH4.3}
%%-----------------------%%

Finally we show the estimate on $\widehat{\Pi}_x^{\epsilon,h;p}(\tau)$ for any $\tau\in\dot{\mcB}_{i+1}\setminus\dot{\mcB}_i$ and $p\in(2,\infty]$. We saw above that the order of the Taylor expansion in the definition of $\widehat{\Pi}_x^{\epsilon,h;p}(\tau)$ varies with $p$; this difference is described by the identity \eqref{BH:eq:L2toLp} below. (A particular case of this identity was proved by Bruned \& Nadeem in their work \cite{BrunedNadeemDiagramFree}, Proposition 3.7 therein.)

\medskip

\begin{lem}\label{lem:BH4.3.1}
\cite[Lemma 6 \& Lemma 16]{RandomModel}
For any $\tau\in\dot{\mcB}_{i+1}$ there exist some finite subsets $\{\sigma_j\}\subset\dot{\mcB}_i ,\, \{\eta_j\}\subset\mcB_i$ and $\{\lambda_j\}\subset\bbR$ satisfying the following properties.
\begin{itemize} \setlength{\itemsep}{0.1cm}
\item[--]
One has $r_p(\tau)=r_p(\sigma_j)+r_\infty(\eta_j)$ and $r_\infty(\sigma_j)\le0<r_2(\sigma_j)$ for each $j$.
\item[--]
For any $\epsilon$ and $h$ there exist some functions $\big\{ x\mapsto f_x^{\epsilon,h}(\sigma_j) \big\}$ such that, for any $p\in[2,\infty]$, we have the identity
\begin{align} \label{BH:eq:L2toLp}
\widehat{\Pi}_x^{\epsilon,h;p}(\tau) = \widehat{\Pi}_x^{\epsilon,h;2}(\tau) + \sum_j \lambda_j \, {\bf1}_{\{r_{p}(\sigma_j)\le0\}} \, f_x^{\epsilon,h}(\sigma_j) \, \overline{\Pi}_x^{\epsilon}(\eta_j).
\end{align}
\item[--]
Denote by $p(\sigma_j)$ the unique $q\in(2,\infty]$ satisfying $r_q(\sigma_j)=0$.
There exists a polynomial $P$ such that for any $h\in H$ and any compact set $C\subset\bbR^d$ we have
$$
\|f_x^{\epsilon,h}(\sigma_j)\|_{L_x^{p(\sigma_j)}(C)}
\le P\big(\|\widehat{\Pi}^{\epsilon,h;p(\sigma_j)}\|_{\mcB_i\cup\dot{\mcB}_i;\overline{C}}\big).
$$
\end{itemize}
\end{lem}

\medskip

The above lemma actually holds for some arbitrary preparation maps, not only for maps of the specific form \eqref{EqTildeR}. Therefore we do not have to go through the extended $\mathfrak{o}$-decoration as in \cite{HS}.

Then we can roughly check that the desired estimate of $\widehat{\Pi}_x^{\epsilon,h;p}(\tau)$ follows from \eqref{BH:eq:L2toLp}. From the assumption $\text{\bf bd}(\dot{\mcB}_{i+1},2)$ the first term of the right hand side of \eqref{BH:eq:L2toLp} satisfies the $B_{2,\infty}^{r_{2}(\tau)}$-type estimate
$$
\big\|(\widehat{\Pi}_x^{\epsilon,h;2}\tau)(\varphi_x^\lambda)\big\|_{L_x^2}\lesssim \lambda^{r_2(\tau)}.
$$
For the remaining terms, the assumption $\text{\bf bd}(\mcB_i)$ and Lemma \ref{lem:BH4.3.1} imply the $B_{p(\sigma_j),\infty}^{r_{\infty}(\eta_j)}$-type estimates
$$
\big\|f_x^{\epsilon,h}(\sigma_j) \big(\overline{\Pi}_x^{\epsilon}\eta_j\big)(\varphi_x^\lambda)\big\|_{L_x^{p(\sigma_j)}}\lesssim \lambda^{r_\infty(\eta_j)}.
$$
The space $B_{2,\infty}^{r_{2}(\tau)}$ is embedded into $B_{p,\infty}^{r_{2}(\tau)-d(1/2-1/p)}=B_{p,\infty}^{r_{p}(\tau)}$. Moreover the space $B_{p(\sigma_j),\infty}^{r_{\infty}(\eta_j)}$ is also embedded into $B_{p,\infty}^{r_{p}(\tau)}$ because
\begin{align*}
r_{\infty}(\eta_j)-d\bigg(\frac1{p(\sigma_j)}-\frac1p\bigg)
&=r_{p}(\tau)-r_{p}(\sigma_j)-d\bigg(\frac1{p(\sigma_j)}-\frac1p\bigg)\\
&=r_{p}(\tau)-r_{p(\sigma_j)}(\sigma_j)
=r_{p}(\tau).
\end{align*}
These formal argument suggest that $\widehat{\Pi}_x^{\epsilon,h;p}(\tau)$ should satisfy a $B_{p,\infty}^{r_{p}(\tau)}$-type estimate. The following statement shows that this is indeed the case.

\medskip

\begin{lem}\label{lem:BH4.3.3}
\cite[Lemma 17]{RandomModel}
For any $p\in[2,\infty]$, there exist a polynomial $P$ and a finite subset $I_p\subset[2,\infty]$ such that for any compact set $C\subset\bbR^d$ we have
$$
\|\widehat{\Pi}^{\epsilon,h;p}\|_{\dot{\mcB}_{i+1};C}
\le P\bigg(\max_{q\in I_{p}}\|\widehat{\Pi}^{\epsilon,h;q}\|_{\mcB_i\cup\dot{\mcB}_i;\overline{C}}\bigg).
$$
\end{lem}

\medskip

Lemmas \ref{lem:BH4.1.2}, \ref{lem:BH4.3.1} and \ref{lem:BH4.3.3} are stated somewhat inaccurately. The values of $r_p$ on both sides of the inequalities described in these lemmas must be slightly different. Precisely, if $r_p$ on the right hand side is defined for a fixed $\alpha_0$, then $r_p$ on the left hand side is redefined for any smaller choice of $\alpha_0$. In \cite{RandomModel} such an adjustment is performed by introducing an additional parameter -- denoted $\varepsilon>0$ therein while the role of $\epsilon$ here is played by an integer parameter $n$ in \cite{RandomModel}. For the sake of simplicity we omit such technical details here.

\end{document}